\documentclass[english]{article}
\usepackage{amssymb}
\usepackage{amsmath}
\usepackage{babel}

\def\thebibliography#1{\section*{References}\list
  {[\arabic{enumi}]}{\settowidth\labelwidth{[#1]}\leftmargin\labelwidth
    \advance\leftmargin\labelsep
    \usecounter{enumi}}
    \def\newblock{\hskip .11em plus .33em minus -.07em}
    \sloppy
    \sfcode`\.=1000\relax}

\newcommand{\refbook}[3]{{\sc #1}{\em\ #2}{\ #3}}
\newcommand{\refer}[5]{{\sc #1}{\ #2}{\em\ #3}{\bf\ #4}{\ #5}}

\newtheorem{lem}{Lemma}[section]
\newtheorem{cor}[lem]{Corollary}
\newtheorem{teo}[lem]{Theorem}
\newtheorem{os}[lem]{Remark}
\newtheorem{defi}[lem]{Definition}
\newtheorem{prop}[lem]{Proposition}

\newcommand{\qed}{\thinspace\null\nobreak\hfill\hbox{\vbox{\kern-.2pt\hrule
 height.2pt depth.2pt\kern-.2pt\kern-.2pt \hbox to2.5mm{\kern-.2pt\vrule
 width.4pt \kern-.2pt\raise2.5mm\vbox to.2pt{}\lower0pt\vtop
 to.2pt{}\hfil\kern-.2pt \vrule
 width.4pt \kern-.2pt}\kern-.2pt\kern-.2pt\hrule height.2pt depth.2pt
 \kern-.2pt}}\par\medbreak}

\newcommand{\R}{\mathbb{R}}

\newcommand{\C}{\mathbb{C}}
\newcommand{\Tt}{\big( T(t) \big)_{t\geq 0}}
\newcommand{\N}{\mathbb{N}}
\newcommand{\Z}{\mathbb{Z}}

\newcommand{\eps}{\varepsilon}

\newcommand{\ov}{\overline}

\newcommand{\ds}{\displaystyle}
\newcommand{\supp}{\emph{supp\,}}

\date{}
\begin{document}

\title{Weighted Calder\'on-Zygmund and Rellich inequalities in  $L^p$}
\author{G. Metafune \thanks{Dipartimento di Matematica ``Ennio De
Giorgi'', Universit\`a del Salento, C.P.193, 73100, Lecce, Italy.
e-mail: giorgio.metafune@unisalento.it, chiara.spina@unisalento.it} \qquad M. Sobajima \thanks{Department of Mathematics, Tokyo University of Science, Japan and Dipartimento di Matematica ``Ennio De
Giorgi'', Universit\`a del Salento, C.P.193, 73100, Lecce, Italy. email: motohiro.sobajima@unisalento.it} \qquad C. Spina \footnotemark[1]}

\maketitle
\begin{abstract}
We find necessary and sufficient conditions for the validity of weighted Rellich and Calder\'on-Zygmund inequalities in $L^p$, $1 \le p \le \infty$, in the whole space and in the half-space with Dirichlet boundary conditions. General operators like $L=\Delta+c\frac{x}{|x|^2}\cdot\nabla-\frac{b}{|x|^2}$ are considered. We compute best constants in some situations.

\bigskip\noindent
Mathematics subject classification (2010): 26D10, 35PXX, 47F05.
\par

\noindent Keywords: Rellich inequalities, Calder\'on-Zygmund inequalities, spectral theory.
\end{abstract}

\section{Introduction}
 In 1956, Rellich proved the inequalities 
$$\left(\frac{N(N-4)}{4}\right)^2\int_{\R^N}|x|^{-4}|u|^2\, dx\leq \int_{\R^N}|\Delta u|^2\, dx$$
for $N\not =2$ and for every $u\in C_c^\infty (\R^N\setminus\{0\})$, see \cite{rellich}.
These inequalities have been then extended to $L^p$-norms:  in 1996, Okazawa proved the validity of
$$\left(\frac{N}{p}-2\right)^p\left(\frac{N}{p'}\right)^p\int_{\R^N}|x|^{-2p}|u|^p\, dx\leq \int_{\R^N}|\Delta u|^p\, dx$$ for $1<p<\frac{N}{2}$ (see \cite{okazawa} and also \cite{Kova-Pere-Seme}) showing also the optimality of the constants. \\

Weighted Rellich inequalities have also been studied. In 1998, Davies and Hinz (\cite[Theorem 12]{davi-hinz}) obtained for $N\geq 3$ and for $2-\frac{N}{p}<\alpha<2-\frac{2}{p}$ 
\begin{equation} \label{WRI}
C(N,p,\alpha)\int_{\R^N}|x|^{(\alpha-2)p}|u|^p\, dx\leq \int_{\R^N}|x|^{\alpha p}|\Delta u|^p\, dx
\end{equation}
with the optimal constants $C(N,p,\alpha)=\left(\frac{N}{p}-2+\alpha\right)^p\left(\frac{N}{p'}-\alpha\right)^p$.
Later Mitidieri showed that  (\ref{WRI}) holds in the wider range $2-\frac{N}{p}<\alpha<N-\frac{N}{p}$ and with the same constants, see  \cite[Theorem 3.1]{mitidieri}.\\
In recent papers Ghoussoub and Moradifam and  Caldiroli and Musina, see \cite{gh-mo}, \cite{caldiroli}, improved  weighted Rellich inequalities for $p=2$ by giving necessary and sufficient conditions on $\alpha$ for the validity of (\ref{WRI}) and finding also the optimal constants $C(N,2,\alpha)$. In particular in \cite{caldiroli} it is  proved that  (\ref{WRI}) is verified for $p=2$ if and only if $\alpha\neq N/2+n$, $\alpha \neq -N/2+2-n$ for every $n\in\N_0$. This approch makes use of the so called Emden-Fowler transform which reduces the operator $|x|^\alpha \Delta$ in $\R^N$ to a uniforly elliptic operator in the cylinder $\R \times S^{N-1}$ and Rellich inequalities to spectral inequalities for the Laplace Beltrami $\Delta_0$ on $S_{N-1}$. We also refer to \cite[Section 3]{gh-mo} where results similar to \cite{caldiroli} have been obtained under the restriction $\alpha \ge (4-N)/2$ and with different methods. \\

In this paper we extend the results in \cite{caldiroli}, \cite{gh-mo} to $1 \le p \le \infty$, computing also best constants in some cases. We show that (\ref{WRI}) holds if and only if $\alpha \neq N/p'+n$, $\alpha \neq -N/p+2-n$ for every $n \in \N_0$.
 Moreover, we use Rellich inequalities to find necessary and sufficient conditions for the validity of weighted Calder\'on-Zygmund estimates when $1<p<\infty$

 \begin{equation} \label{WCZI}
\int_{\R^N}|x|^{\alpha p}|D^2u|^p\, dx\leq C\int_{\R^N}|x|^{\alpha p}|\Delta u|^p\, dx
\end{equation}
for $u \in C_c^\infty (\R^N \setminus \{0\})$.

Weighted Calder\'on-Zygmund inequalities are well-known in the literature, in the framework of singular integrals. In 1957 Stein (see \cite{stein}) proved the inequalities
\begin{equation} \label{Stein}
\||x|^\alpha Tf\|_p \le C\||x|^\alpha f\|_p
\end{equation}
for $1<p<\infty$, $-N/p < \alpha <N/p'$, where $T$ is the Calder\'on-Zygmund kernel corresponding to the operator $D^2\Delta^{-1}$.
Subsequent generalizations of the above result can be found in the papers of Kree, Muckenhoupt and Wheeden (see \cite{kree}, \cite{muck-wheed}) where more general kernels are treated.
Taking  $u\in C_c^\infty(\R^N \setminus \{0\})$ and setting $f=\Delta u$, inequalities (\ref{Stein}) imply that
\begin{equation*} 
\||x|^\alpha D^2 u\|_p \leq C \||x|^{\alpha} \Delta u\|_p.
\end{equation*}
 However the last inequalities can hold also when (\ref{Stein}) fail, that is outside of the range $-N/p < \alpha <N/p'$, since $f$ has compact support whenever $u$ has but the converse is clearly false. In particular, the condition $\alpha >-N/p$ is needed for the integrability of $|x|^\alpha Tf$ near the origin, whereas $\alpha<N/p'$ is needed for the integrability at infinity, if $Tf$ behaves like $|x|^{-N}$.
 We find that (\ref{WCZI}) holds if and only if $\alpha\neq N/p'+n$ for every $n \in \N_0$ and , $\alpha \neq -N/p+2-n$ for every $n\in\N, n \ge 2$.

We consider also more general operators 
$$L=\Delta+c\frac{x}{|x|^2}\cdot\nabla-\frac{b}{|x|^2}$$ with $b, c \in \C$  and investigate the validity of weighted Rellich inequalities of the form
\begin{equation} \label{WRI-L}
C(N,p,\alpha,b,c)\int_{\R^N}|x|^{(\alpha-2)p}|u|^p\, dx \leq \int_{\R^N}|x|^{\alpha p}|L u|^p\, dx
\end{equation}
for $u \in C_c^\infty (\R^N \setminus \{0\})$ and $1\leq p < \infty$.
We prove necessary and sufficient conditions on $\alpha$ for the validity of (\ref{WRI-L}) and, in certain cases, we explicitely compute the best constants. 
\bigskip

Let us describe more analitically the content of the paper. In Section 2 we prove Rellich inequalities for $p=2$ by using orthogonal decomposition in spherical harmonics and Hardy inequalities. In particular we recover with a different method the results of Caldiroli and Musina quoted in the Introduction. This approach had the advantage of being elementary but leads to the unnecessary condition (\ref{condition-prop3.1}) which will be removed in Subsection 2.3 using the  spectral arguments which will be the basis for the $L^p$-analysis. 

The $L^p$ case is treated in Section 3, including the endpoints $p=1,\infty$. With the change of unknown $v=|x|^{\alpha-2}u$, Rellich inequality (\ref{WRI-L}) is equivalent to the inequality
$$
\|\tilde{L}v-bv\|_p \geq C\|v\|_p
$$
where
$$
\tilde L=|x|^2\Delta +(4-2\alpha+c)x\cdot\nabla +(2-\alpha)(N-\alpha+c)
$$
(with the same constant $C$) which holds if and only if $b$ does not belong to the spectrum of ${\tilde L}$.
With this approach Rellich inequalities are then reduced to  spectral problems for singular operators which are analyzed in detail in Section 5, Section 6  using tensor products arguments. Note that an operator like $A=|x|^2\Delta+cx\cdot \nabla$ can be written in spherical coordinates as 
$$A=\rho^2\frac{\partial^2}{\partial\rho^2}+(N-1+c)\rho \frac{\partial }{\partial \rho}+\Delta_0,$$
$\Delta_0$ being the Laplace-Beltrami on unit sphere, and turns out to be the sum of two (commuting) operators acting on independent variables.

Section 4  is devoted to the analysis of weighted Calder\'on-Zygmund estimates. We show that, apart from special values of $\alpha,p$, Rellich and Calder\'on-Zygmund inequalities are equivalent and find all values of $\alpha$ for which (\ref{WCZI}) hold.

In Section 7 we present special cases of our results and generalizations. In particular we show that Rellich and Calder\'on-Zygmund inequalities can hold on subspaces of $L^p(\R^N)$ defined by sets of spherical harmonics, even though they fail in the whole space. For example, the classical Rellich inequalities (that is with $\alpha=b=c=0$) which fail for $p=N/2,N$ continue to hold when $p=N/2$ for functions having zero mean on $S_{N-1}$ and when $p=N$ for functions orthogonal to spherical harmonics of order 1. Similar remarks hold for $p=1,\infty$: in $L^1$ they hold for functions with zero mean as above and in $L^\infty$ for functions orthogonal to spherical harmonics of degree 2. Moreover we find the best constants on special subspaces defined by shperical harmonics of a fixed order and give estimates from above and from below in the case of the whole space.

The case of the half-plane is analyzed in detail showing that Rellich and Calder\'on-Zygmund inequalities hold for functions vanishing at the boundary in some cases where they fail in the whole space.

In Appendix A we collect and prove the Hardy type inequalities we need in the paper, showing also the optimality of the constants. In Appendix B we briefly analyze the singular operator $A=|x|^2\Delta +cx\cdot \nabla$ in spaces of continuous functions to clarify the nature of the singularities $0, \infty$ from the point of view of the underlying stochastic process.
In Appendix C we recall and proof a result on the norm of the tensor product of two operators.

\bigskip \noindent
{\bf Acknowledgment.} The authors thank Roberta Musina for many comments and suggestions on the paper.

\bigskip
\noindent\textbf{Notation.}
We denote by $\N_0=\N\cup\{0\}$ the natural numbers including 0. $C_b(\Omega)$ is the Banach space of all continuous and bounded functions in $\Omega$, endowed with the sup-norm, $C_0(\Omega)$ its subspace consiting of functions vanishing at the boundary of $\Omega$. $C_c^\infty(\Omega)$ denotes the space of infinitely continuously differentiable functions with compact support in $\Omega$. $C_0^0(\R^N)$ stands for the Banach space of all continuous functions in $\R^N$ vanishing at $0,\infty$.
The unit sphere $\{\|x\|=1\}$ in $\R^N$ is denoted by $S_{N-1}$ and  $\Sigma$ is an open $C^2$-set of $S_{N-1}$; $\Delta_0$ is the Laplace-Beltrami operator endowed with Dirichlet boundary conditions if $\Sigma\neq S_{N-1}$. We adopt standard notation for $L^p$ and Sobolev spaces . When $p=\infty$ we write $L^\infty (\Omega)$ for $C_0(\Omega)$.

\section{ Rellich inequalities in $L^2(\R^N)$}
\subsection{The Laplacian  in spherical coordinates} \label{preliminaries}
We introduce spherical coordinates
\begin{equation} \label{spherical}
\left\{
\begin{array}{ll}
  x_1=\rho \cos\theta_1\sin \theta_2\ldots\sin \theta_{N-1}\\ 
  x_2=\rho \sin\theta_1\sin \theta_2\ldots\sin \theta_{N-1}\\
 \vdots\\
 x_n=\rho\cos \theta_{N-1}\end{array}\right.
\end{equation}
where $\theta_2,\ldots,\theta_{N-1}$ range from $0$ to $\pi$ and $\theta_1$ ranges from $0$ to $2\pi$.
The Laplace operator is then  given by 
$$\Delta= \frac{\partial^2}{\partial\rho^2}+\frac{N-1}{\rho}\frac{\partial }{\partial \rho}+\frac{1}{\rho^2}\Delta_0$$
where 
$$\Delta_0=\frac{1}{\sin^{N-2}\theta_{N-1}}\frac{\partial}{\partial\theta_{N-1}}\sin^{N-2}\theta_{N-1}\frac{\partial}{\partial\theta_{N-1}}+\ldots+\frac{1}{\sin^2\theta_{N-1}\cdots\sin^2\theta_2}\frac{\partial^2}{\partial\theta_1^2}$$
 is the Laplace-Beltrami on the unit sphere $S_{N-1}$, see \cite[Chapter IX]{vilenkin}) . 
If $u(x)=u(\rho,\omega)\in C_c^\infty(\R^N \setminus \{0\})$, $\rho\in[0,\infty[$, $\omega\in S_{N-1}$,
by \cite[Ch. 4, Lemma 2.18]{stein-weiss}, 
$$u(x)=\sum_{n=0}^\infty c_n(\rho)P_n(\omega)$$ 
in $L^2(\R^N)$, where $(P_n)$ is a complete orthonormal  system of spherical harmonics and $$c_n(\rho)=\int_{S_{N-1}}u(\rho,\omega)P_n(\omega)d\omega.$$

 By the regularity of $u$ it follows that $c_n(\rho)\in C_c^\infty(]0,\infty[)$.
We recall that a spherical harmonic of order $n$ is the restriction to $S_{N-1}$ of a homogenuous harmonic polynomial of degree $n$.

\begin{lem} 
 Let $P_n$ be a spherical harmonics of order $n$ on $S_{N-1}$. Then for every $n\in\N_0$
$$\Delta_0 P_n=-(n^2+(N-2)n)P_n.$$ 
The values $\lambda_n=n^2+(N-2)n$ are the eigenvalues of the Laplace-Beltrami operator $-\Delta_0$ on $S_{N-1}$. The corresponding eigenspace consists of all spherical harmonics of order $n$ and has dimension $d_n$ where $d_0=1, d_1=N$ and
$$
d_n= \binom{N+n-1}{n}-\binom{N+n-3}{n-2}.
$$
for $n \ge 2$.
\end{lem}
It follows that, if $u\in C_c^\infty(\R^N \setminus \{0\})$, $u(x)=\sum_{n=0}^\infty c_n(\rho)P_n(\omega)$, then
\begin{equation} \label{represent}
\Delta u(\rho,\omega)=\sum_{n=0}^\infty \left(c''_n(\rho)+\frac{N-1}{\rho}c'_n(\rho)-\lambda_n\frac{c_n(\rho)}{\rho^2}\right)P_n(\omega), 
\end{equation}
where the eigenvalues $\lambda_n$ are repeated according to their multiplicity.

\subsection{Rellich inequalities in $L^2$: Part I}
In this section we prove weighted Rellich inequalities 
for general operators
of the form 
\begin{equation} \label{defL}
Lu:=\Delta u+c\frac{x}{|x|^2}\cdot\nabla u-\frac{b}{|x|^{2}}u,
\end{equation}
where $c ,b\in \R$. 
The proof is based on integration by parts and Hardy's inequalities but leads to condition (\ref{condition-prop3.1}) which will be removed in Part II.
In order to shorten the notation we set 
\begin{equation} \label{defgamma2}
\gamma_2(\alpha,c)=\Bigl(\frac{N}{2}-2+\alpha\Bigr)\Bigl(\frac{N}{2}-\alpha+c\Bigr)=\Bigl(\frac{N}{2}-1+\frac{c}{2}\Bigr)^2-\Bigl(1-\alpha+\frac{c}{2}\Bigr)^2.
\end{equation}

\begin{prop} \label{Rellich}
Let $N\geq 2$,  $\alpha\in \R$ such that 
$b+\lambda_n+\gamma_2(\alpha,c)\neq 0$ for every $n \in \N_0$. 
If
\begin{equation}\label{condition-prop3.1}
b+\gamma_{2}(\alpha,c)+2\Bigl(1-\alpha+\frac{c}{2}\Bigr)^2=b+\Bigl(\frac{N}{2}-1+\frac{c}{2}\Bigr)^2+\Bigl(1-\alpha+\frac{c}{2}\Bigr)^2\geq0,
\end{equation}
then for every $u\in C_c^\infty(\R^N \setminus \{0\})$,
$$\int_{\R^N}|x|^{2\alpha}|L u|^2\, dx\geq C^2(N,\alpha,b,c)\int_{\R^N}|x|^{2\alpha-4}|u|^2\,dx$$
 where  
$$C^2(N,\alpha,b,c)=\min_{n\in\N_0}\left (b+\lambda_n+\gamma_2(\alpha,c)\right )^2>0.$$
\end{prop}
{\sc Proof.}
Let $u\in C_c^\infty(\R^N\setminus\{0\})$ 
and set $v=|x|^{\alpha+\frac{N}{2}-2}u$. 
Then we note that 
\[
|x|^{\alpha}Lu
=
|x|^{-\frac{N}{2}+2}
\Delta v
+
(4-2\alpha-N+c)|x|^{-\frac{N}{2}}x\cdot\nabla v
-
\left[
b+
\gamma_{2}(\alpha,c)
\right]|x|^{-\frac{N}{2}}v.
\]
Expanding $v$ in  spherical harmonics 
\begin{equation*}
v(x)=v(\rho,\omega)=\sum_{n=0}^{\infty}d_n(\rho)P_n(\omega), 
\end{equation*}
by (\ref{represent}) we have 
\begin{equation}  \label{laplacian}
\int_{\R^N}
  |x|^{2\alpha}|Lu|^2
\,dx
=
\sum_{n=0}^\infty
\int_{0}^\infty
  \rho^3\left|d_n''(\rho) +\frac{\tilde{c}}{\rho} d_n'(\rho) -\frac{(\tilde{b}+\lambda_n)}{\rho^2}d_n(\rho)\right|^2
\,d\rho, 
\end{equation}
and
\begin{equation}  \label{norm-u}
\int_{\R^N}|u|^2|x|^{2\alpha-4}\,dx=
\sum_{n=0}^\infty \int_0^\infty \rho^{-1} |d_n(\rho)|^2\,d\rho,
\end{equation}
where 
\begin{equation}\label{coeff_L2}
\tilde{c}:=3-2\alpha+c, 
\quad 
\tilde{b}:=
b+\gamma_2(\alpha,c) 
\end{equation}
and the eigenvalues $\lambda_n$ are repeated according to their multiplicity.
Recalling that $d_n\in C_c^\infty((0,\infty))$ and integrating by parts we obtain
\begin{eqnarray*}
& &
\int_{0}^\infty
  \rho^3\left|d_n''(\rho) +\frac{\tilde{c}}{\rho} d_n'(\rho) -\frac{(\tilde{b}+\lambda_n)}{\rho^2}d_n(\rho)\right|^2
\,d\rho
\\&=&
\int_0^\infty \rho^{3}|d''_n(\rho)|^2\,d\rho
+
[\tilde{c}^2-2\tilde{c}+2(\tilde{b}+\lambda_n)]
\int_0^\infty \rho
|d'_n(\rho)|^2\,d\rho
\\
&&+
(\tilde{b}+\lambda_n)^2 
\int_0^\infty \rho^{-1}
|d_n(\rho)|^2\,d\rho.
\end{eqnarray*}
To estimate the second order derivative 
we use the following Hardy inequality 
\[
\int_0^\infty \rho^3|w'(\rho)|^2\,d\rho
\geq 
\int_0^\infty \rho|w(\rho)|^2\,d\rho, \quad w\in C_c^\infty((0,+\infty))
\]
(see Proposition \ref{hardy}). 
Taking $w=d'$ implies that
\begin{align}
\nonumber 
&\int_{0}^\infty
  \rho^3\left|d_n''(\rho) +\frac{\tilde{c}}{\rho} d_n'(\rho) -\frac{(\tilde{b}+\lambda_n)}{\rho^2}d_n(\rho)\right|^2
\,d\rho
\\& \nonumber\geq 
\left[(\tilde{c}-1)^2+2(\tilde{b}+\lambda_n)\right]\int_0^\infty \rho|d'_n(\rho)|^2\,d\rho
+
(\tilde{b}+\lambda_n)^2\int_0^\infty\rho^{-1}|d_n(\rho)|^2\,d\rho.
\end{align}
By virtue of \eqref{condition-prop3.1}, for every $n\in\N_0$, the coefficient of the first term on the right-hand side 
of the above estimate is nonnegative
\[
(\tilde{c}-1)^2+2(\tilde{b}+\lambda_n)
=
2\left[b+\gamma_{2}(\alpha,c)+2\Bigl(1-\alpha+\frac{c}{2}\Bigr)^2\right]
+
2\lambda_n
\geq 0.
\]
Thus noting that $\tilde{b}+\lambda_n=b+\lambda_n+\gamma_{2}(\alpha,c)$, 
we have 
\[
\int_{0}^\infty
  \rho^3\left|d_n''(\rho) +\frac{\tilde{c}}{\rho} d_n'(\rho) -\frac{(\tilde{b}+\lambda_n)}{\rho^2}d_n(\rho)\right|^2
\,d\rho
\geq 
[b+\lambda_n+\gamma_2(\alpha,c)]^2
\int_0^\infty\rho^{-1}|d_n(\rho)|^2\,d\rho.
\]
Consequently, we obtain 
\begin{align*}  
\int_{\R^N}|x|^{2\alpha}|L u|^2\,dx
&=\sum_{n=0}^\infty \int_{0}^\infty
  \rho^3\left|d_n''(\rho) +\frac{\tilde{c}}{\rho} d_n'(\rho) -\frac{(\tilde{b}+\lambda_n)}{\rho^2}d_n(\rho)\right|^2
\,d\rho
\\&\geq
\sum_{n=0}^\infty \left[(b+\lambda_n+\gamma_2(\alpha,c))^2\int_0^\infty\rho^{-1}|d_n(\rho)|^2\,d\rho\right]
\\
&\geq
\min_{n\in\N_0}\left[b+\lambda_n+\gamma_2(\alpha,c)\right]^2\sum_{n=0}^\infty  \int_0^\infty\rho^{-1}|d_n(\rho)|^2\,d\rho\\
&\geq
\min_{n\in\N_0}\left[b+\lambda_n+\gamma_2(\alpha,c)\right]^2
\int_{\R^N}|x|^{2\alpha-4}|u|^2\,dx.
\end{align*}
\qed

\begin{os}
{\rm The optimality of the constant $C(N,\alpha,b,c)$ 
 will be proved in the next section.}
\end{os}

We point out that, in correspondence of $b=c=0$, Proposition \ref{Rellich} provides, with an alternative proof, the same result contained in  \cite[Theorem 4.1]{caldiroli}, see also \cite[Theorem 3.14]{gh-mo}.

\begin{cor} \label{CM2} If $b=c=0$, that is if $L=\Delta$, then Rellich  inequalities hold in $L^2(\R^N)$ if and only if for every $n \in \N_0$
$$
\alpha \neq \frac{N}{2}+n, \qquad \alpha \neq -\frac{N}{2}+2-n.
$$
Moreover, the best constant is given by
$$
C=\min_{n \in \N_0} \left |\left (n+\frac{N}{2} -1\right )^2-(1-\alpha)^2 \right |.
$$
\end{cor}
{\sc Proof. }
Indeed  condition (\ref{condition-prop3.1}) is satisfied since
$\gamma_2(\alpha,0)+2(1-\alpha)^2= \left(\frac{N}{2}-1\right)^2+(1-\alpha)^2\geq 0$.  Next observe that $\lambda_n=n^2+(N-2)n$ yields 
$$\lambda_n+\gamma_2(\alpha,0)=\left (n+\frac{N}{2} -1\right )^2-(1-\alpha)^2=\left (\frac{N}{2}-2+\alpha+n\right )\left (\frac{N}{2}-\alpha+n\right ),$$ hence  the requirement $\lambda_n+\gamma_2(\alpha,0)\neq 0$ leads to the statement.
\qed

\subsection{Rellich inequalities in $L^2$: Part II}

In this section we prove Rellich inequalities for operators as in (\ref{defL}) using spectral arguments which will be the basis for the $L^p$ analysis. In particular we remove condition (\ref{condition-prop3.1}) and compute best constants.
Here we restrict ourselves to the case where $b$ and $c$ are real. 

To state the result in this section, we introduce
\begin{equation} \label{P2}
{\cal P}_{2,\alpha,c}:=
\left\{\lambda=-\xi^2+i\xi(2-2\alpha+c)-\gamma_2(\alpha,c)\;;\;\xi\in \R\right\},
\end{equation}
where $\gamma_2(\alpha,c)$ is defined in (\ref{defgamma2}).

\begin{teo}\label{WRellichL2}
Let $N\geq 2$ and $\alpha, b,c\in \R$. Then Rellich inequalities
\begin{equation}
\label{WRellichL2-ineq}
\int_{\R^N}|x|^{2\alpha}|L u|^2\,dx
\geq 
C^2
\int_{\R^N}|x|^{2\alpha-4}|u|^2\,dx
\end{equation}
hold for every $u\in C_c^\infty(\R^N\setminus\{0\})$ and with $C>0$ independent of $u$, if and only if $b+\lambda_n \not\in {\cal P}_{2,\alpha,c}$ for every $n \in \N_0$. 
In such a case the optimal constant $C$ is given by   
$C^2(N,\alpha,b,c):=\min_{n\in \N_0} C^2_n(N,\alpha,b,c)>0$ where 
\[
C^2_n(N,\alpha,b,c)
:=
\begin{cases}
\Bigl(
b+\lambda_n+\gamma_2(\alpha,c)
\Bigr)^2
\\[8pt]
\hspace{130pt}{\rm if}\ 
b+\lambda_n+\Bigl(\frac{N}{2}-1+\frac{c}{2}\Bigr)^2+\Bigl(1-\alpha+\frac{c}{2}\Bigr)^2 \ge 0, 
\\[8pt]
4\Bigl(1-\alpha+\frac{c}{2}\Bigr)^2
\left[-b-\Bigl(\frac{N}{2}-1+\frac{c}{2}\Bigr)^2-\lambda_n\right]
\\[8pt]
\hspace{130pt}{\rm if}\ 
b+\lambda_n+\Bigl(\frac{N}{2}-1+\frac{c}{2}\Bigr)^2+\Bigl(1-\alpha+\frac{c}{2}\Bigr)^2<0.
\end{cases}
\]
\end{teo}

\begin{os} {\rm The condition  $b+\lambda_n \not \in {\cal P}_{2,\alpha, c}$ for every $n\in \N_0$, can be written in a simpler form since $b \in \R$. In fact, if $2-2\alpha-c\not = 0$, then ${\cal P}_{2,\alpha,c}$ is a non degenerate parabola with vertex at $(-\gamma_2(\alpha,c),0)$ and the above condition reads $b+\gamma_2(\alpha,c)+\lambda_n \not = 0$ for every $n\in N_0$, as in the statement of Proposition \ref{Rellich}. However, if $2-2\alpha-c= 0$, then ${\cal P}_{2,\alpha,c}$  coincides with the semiaxis $]-\infty,-\gamma_2(\alpha,c)]$ and the condition becomes $b+\gamma_2(\alpha,c)>0$ (recall that $\lambda_n \ge 0$ and $\lambda_0=0$).}
\end{os}

Recalling that 
$$
\gamma_2(\alpha,c)+2\left (1-\alpha+\frac{c}{2}\right )^2=\Bigl(\frac{N}{2}-1+\frac{c}{2}\Bigr)^2+\Bigl(1-\alpha+\frac{c}{2}\Bigr)^2,
$$
we note that the value in condition (\ref{condition-prop3.1}) is the same where the above formula for $C_n(N,\alpha,b,c)$  changes shape.  Before proving Theorem \ref{WRellichL2} we state the following elementary lemma.
\begin{lem}\label{dist}
Let $\kappa\in \R$ and 
${\cal P}$ be the parabola  
\[
{\cal P}:= \{-\xi^2+2i \kappa\xi\;;\;\xi\in \R\}. 
\]
Then for every $\lambda\in \R$,  
\[
{\rm dist}(\lambda,{\cal P})^2
=
\begin{cases}
\lambda^2&{\rm if}\ \lambda\geq-2\kappa^2, 
\\
4\kappa^2(-\lambda-\kappa^2)&{\rm if}\ \lambda<-2\kappa^2. 
\end{cases}
\]
\end{lem}
Note that the focus of the parabola is in $(-\kappa^2,0)$ and that ${\cal P}=]-\infty,0]$ when $\kappa=0$. 

\smallskip

{\sc Proof of Theorem \ref{WRellichL2}.}
Let $u\in C_c^\infty(\R^N\setminus\{0\})$. 
Noting that 
\eqref{laplacian}, \eqref{norm-u} and \eqref{coeff_L2}
and 
employing change of variables from $\rho$ to $e^s$ 
and putting $w_n(s):=d_n(e^s)\in C_c^\infty(\R)$ imply
\begin{align*}
\int_{\R^N}
  |x|^{2\alpha}|Lu|^2
\,dx
&=
\sum_{n=0}^\infty
\int_{-\infty}^\infty
  \bigl|w_n'' +(\tilde{c}-1)w_n' -(\tilde{b}+\lambda_n)w_n\bigr|^2
\,ds, 
\\
\int_{\R^N}
  |x|^{2\alpha-4} |u|^2
\,dx
&=
\sum_{n=0}^\infty
\int_{-\infty}^\infty
  |w_n|^2\,ds. 
\end{align*}
Thus \eqref{WRellichL2-ineq} is translated into the following families of spectral  inequalities 
\begin{equation} \label{1d}
\int_{-\infty}^\infty
  \bigl|w_n'' +(\tilde{c}-1)w_n' -(\tilde{b}+\lambda_n)w_n\bigr|^2
\,ds \ge c_n \int_{-\infty}^\infty
  |w_n|^2\,ds,
\end{equation}
$c_n>0$, for the one-dimensional operator  
$\Gamma_0:=d^2/dx^2+(\tilde{c}-1)d/dx$ in $L^2(\R)$. 
Since, using the Fourier transform,  
\begin{align*}
\sigma(\Gamma_0)
&=
\left\{\lambda=-\xi^2+i\xi(\tilde{c}-1)\;;\;\xi\in \R\right\}
\end{align*}
it coincides with its topological boundary, hence every point in the spectrum is in the approximate point spectrum. Since $C_c^\infty (\R)$ is a core for $\Gamma_0$, it follows that   (\ref{1d}) holds  if and only if $\tilde{b}+\lambda_n\not\in \sigma(\Gamma_0)$ (see \cite[Proposition 1.10, Chapter IV]{engel-nagel} for these elementary properties of the approximate point spectrum).
Moreover, since the resolvent of $\Gamma_0$ is a normal operator, the spectral theorem and Lemma \ref{dist} with $\kappa=1-\alpha+c/2 $ give 
\begin{align*}
\|(\Gamma_0 - (\tilde{b}+\lambda_n))^{-1}\|^2
=
{\rm dist} (\tilde{b}+\lambda_n, \sigma(\Gamma_0))^{-2}
=
\left (C_n(N,\alpha,b,c)\right)^{-2},
\end{align*}
that is 
\begin{align*}
\int_{-\infty}^\infty
  |w_n|^2
\,ds
\leq
\frac{1}{C^2_n(N,\alpha,b,c)}
\int_{-\infty}^\infty
  \bigl|
  \Gamma_0 w_n-(\tilde{b}+\lambda_n)w_n\bigr|^2
\,ds.
\end{align*}
Thus we see 
that if 
$\tilde{b}+\lambda_n\not\in \sigma(\Gamma_0)$
for every $n\in \N_0$, then 
\begin{align*}
\int_{\R^N}
  |x|^{2\alpha-4}
  |u|^2
\,dx
&=
\sum_{n=0}^\infty
\int_{-\infty}^\infty
  |w_n|^2
\,ds
\\
&
\leq
\sum_{n=0}^\infty
\left(
\frac{1}{C^2_n(N,\alpha,b,c)}
\int_{-\infty}^\infty
  \bigl|
      \Gamma_0 w_n -(\tilde{b}+\lambda_n)w_n\bigr|^2
\,ds
\right)
\\
&
\leq
\frac{1}{\min_{n\in \N_0}C^2_n(N,\alpha,b,c)}
\int_{\R^N}
  |x|^{2\alpha}|Lu|^2
\,dx.
\end{align*}
This is nothing but the desired inequality. 

\begin{os}
{\rm The constant $C_n(N,\alpha,b,c)$ is optimal for every $n\in \N_0$ and
this implies the optimality of $C(N,\alpha,b,c)$. Actually, $C_n(N, \alpha,b,c)$ is the best constant for which inequalities (\ref{WRellichL2-ineq}) hold when $u=\sum_k c_k(\rho) P_k(\omega)$, where the sum is finite and all the spherical harmonics $P_k$ have order $n$ (hence are eigenfunctions of $-\Delta_0$ with eigenvalue $\lambda_n$). This is easily seen from the proof.}
\end{os}

\bigskip
\begin{os} {\rm In the case of  Schr\"odinger operators, that is when $c=0$, the best constant can be of the form 
$$
4\Bigl(1-\alpha \Bigr)^2
\left[-b-\Bigl(\frac{N}{2}-1\Bigr)^2-\lambda_{n_1}\right]=4\Bigl(1-\alpha \Bigr)^2
\left[-b-\Bigl( n_1+\frac{N}{2}-1\Bigr)^2\right]
$$
 for some  $n_1\in\N_0$.}

\bigskip\noindent
{\rm
Let $b<0$, $n_1\in\N_0$ be such that 
$b+\left(n+\frac{N}{2}-1\right)^2+(1-\alpha)^2<0$ if and only if $n\leq n_1$. 
For $n\leq n_1$, 
$$C^2_n=4\Bigl(1-\alpha\Bigr)^2
\left[-b-\Bigl(n+\frac{N}{2}-1\Bigr)^2\right]\geq C^2_{n_1}$$
and if $n>n_1$, 
$$C^2_n=
\left (b+\Bigl(n+\frac{N}{2}-1\Bigr)^2
\right)^2.$$

We  fix $n_1=0$ and choose $b=-\left(\frac{N}{2}-1\right)^2- 2(1-\alpha)^2-1$. Then $b+\left(n+\frac{N}{2}-1\right)^2+2(1-\alpha)^2\geq 0$ if and only if  $n>0$. 
The inequality
$$C_0^2=4(1-\alpha)^2\left (1+2(1-\alpha)^2\right) \le C_n^2=\left (n^2+n(N-2)-2(1-\alpha)^2-1 \right )^2$$
holds for every  $n >0$ if $N$ is sufficiently large.

}
\end{os}

\section{Rellich inequalities in $L^p(\R^N)$} \label{rellichp}
In this Section we prove Rellich inequalities in $L^p$, $1 \le p \le \infty$. 
As before, we set 
$$Lu=\Delta u+c\frac{x}{|x|^2}\cdot\nabla u-\frac{b}{|x|^2} u.$$ 
The coefficients $c$ and $b$ are allowed to be complex. We shall be able to determine all $\alpha's$ (depending on $N,p,c,b$) for which Rellich inequalities hold but, in contrast with the case $p=2$, we can compute the best constants only under additional conditions. We start with the case where the coefficients $b$ and $c$ are real.

\subsection{Real coefficients}
As in the $L^2$ case we set for $1\le p \le \infty$

\begin{equation} \label{gammap}
\gamma_p(\alpha,c)=\Bigl(\frac{N}{p}-2+\alpha\Bigr)\Bigl(\frac{N}{p'}-\alpha+c\Bigr)=\Bigl(\frac{N}{2}-1+\frac{c}{2}\Bigr)^2-\Bigl(N\Bigl (\frac12-\frac1p\Bigr)+1-\alpha+\frac{c}{2}\Bigr)^2
\end{equation}
and the parabola
\begin{equation} \label{Pp}
{\cal P}_{p,\alpha,c}:=
\left\{\lambda=-\xi^2+i\xi\Bigl (N\Bigl (1-\frac2p \Bigr)+2-2\alpha+c\Bigr)-\gamma_p(\alpha,c)\;;\;\xi\in \R\right\}.
\end{equation}

\begin{teo} \label{Rellich-p}
Let $N\geq 2$, $\alpha, b,\ c\in\R$, $1 \le p \le \infty$.
There exists a positive constant $C=C(N,\alpha,  p, c, b)$ such that
\begin{equation} \label{rellich-op}
\||x|^\alpha Lu\|_p \geq C\||x|^{\alpha-2} |u|\|_p
\end{equation}
holds for every $u\in C_c^\infty(\R^N \setminus \{0\})$ if and only if $b+\lambda_n \not \in {\cal P}_{p,\alpha,c}$ for every $n \in \N_0$.
   If, in addition,  $b+\gamma_p(\alpha,c)>0$
the optimal constant is given by
$C=b+\gamma_p(\alpha,c).$

 \end{teo}

\begin{os} {\rm The condition $b\not\in \bigcup_{n=0}^\infty({\cal P}_{p,\alpha,c}-\lambda_n)$ or $b+\lambda_n \not \in {\cal P}_{p,\alpha, c}$ for every $n\in \N_0$, can be written in a simpler form since $b \in \R$. In fact, if $N(1-2/p) +2-2\alpha+c \not =0$, then ${\cal P}_{p,\alpha,c}$ is a non degenerate parabola with vertex at $(-\gamma_p(\alpha,c),0)$ and the above condition reads $b+\gamma_p(\alpha,c)+\lambda_n \not = 0$ for every $n\in \N_0$ or, equivalently, 
$$
b+\lambda_n+\Bigl(\frac{N}{2}-1+\frac{c}{2}\Bigr)^2 \neq\Bigl(N\Bigl (\frac12-\frac1p\Bigr)+1-\alpha+\frac{c}{2}\Bigr)^2.
$$ However, if $N(1-2/p) +2-2\alpha+c  =  0$, then ${\cal P}_{p,\alpha,c}$  coincides with the semiaxis $]-\infty,-\gamma_p(\alpha,c)]$ and the condition becomes $b+\gamma_p(\alpha,c)>0$.}
\end{os}

{\sc Proof of Theorem \ref{Rellich-p}.} Let $u\in C_c^\infty(\R^N\setminus\{0\})$. Set $v(x)=|x|^{\alpha-2}u(x)$ we observe that
$$|x|^\alpha L u:=\tilde{L}v-bv$$
where 
\begin{equation} \label{deftildeL}
\tilde L=|x|^2\Delta +(4-2\alpha+c)x\cdot\nabla +(2-\alpha)(N-\alpha+c).
\end{equation}
Therefore  (\ref{rellich-op}) is equivalent to the estimate
\begin{equation} \label{disv}
\|\tilde{L}v-bv\|_p \geq C\|v\|_p
\end{equation}
for any $v\in C_c^\infty (\R^N \setminus \{0\})$ or, since $C_c^\infty (\R^N \setminus \{0\})$ is a core for the domain of $\tilde L$, see Proposition \ref{core}, for any $v$ in the domain of $\tilde L$.  Then (\ref{disv}) is true if and only if $b$ does not belong to the  spectrum of  
$\tilde{L}$. Indeed, by Proposition \ref{specRN}
$$\sigma(\tilde{L})=\cup_{n\in\N_0}({\cal P}_{p,\alpha,c}-\lambda_n),$$ where ${\cal P}_{p,\alpha,c}$ is defined in (\ref{Pp}), hence it coincides with its topological boundary. Then $\sigma (\tilde L )$ consists of approximate eigenvalues, that is of all $b$ for which (\ref{disv}) fails, see \cite[Proposition 1.10, Chapter IV]{engel-nagel}.
If $b \not \in \sigma (\tilde L)$, then the optimal constant in (\ref{rellich-op}) is given by 
$$C^{-1}=\|({\tilde L}-b)^{-1}\|_p \ge  (dist(b, \sigma(\tilde{L})))^{-1}$$
hence $C \le dist(b, \sigma(\tilde{L}))$.
Finally, let us assume that  $b+\gamma_p(\alpha,c)>0$. Since by Proposition \ref{dissipativity1}
$$
\|e^{t \tilde L}\|_p \le e^{-t \gamma_p(\alpha,c)}
$$
the resolvent estimate  $\|(\tilde{L}-b)^{-1}\|_p\leq (b+\gamma_{p,\alpha,c})^{-1}=(dist(b, \sigma(\tilde{L})))^{-1}$ also holds and yields $C=b+\gamma_{p,\alpha,c}$. 
\qed 

We specialize the above reult to the case $L=\Delta$, thus obtaining the $L^p$-version of the result of Caldiroli and Musina, \cite{caldiroli}. We note that the extreme point $p=1,\infty$ are allowed.

\begin{teo} \label{CMp} Let $N \ge 2$, $1\le p\le \infty$, $\alpha \in \R$, $b=c=0$. Then Rellich inequalities
$$
\||x|^\alpha \Delta u\|_p \geq C\||x|^{\alpha-2} |u|\|_p
$$
hold for every $u\in C_c^\infty(\R^N \setminus \{0\})$ and for a suitable $C>0$ if and only if 
\begin{equation} \label{condRp}
\alpha \neq \frac{N}{p'}+n, \qquad \alpha \neq -\frac{N}{p}+2-n \qquad {\rm for\  every\ } n \in \N_0.
\end{equation}
Moreover, if $N \ge 3$ and $2-N/p <\alpha <N/p'$, the best constant $C$ is given by 
\begin{equation} \label{bestp}
C=\left (\frac{N}{p}-2+\alpha \right )\left (\frac{N}{p'}-\alpha\right ).
\end{equation}
\end{teo}
{\sc Proof. } The parabola ${\cal P}_{p,\alpha,0}$ degenerates if and only if $\bar{\alpha}=N(1/2-1/p)+1$ and $\gamma_p(\bar{\alpha},0)>0$ if and only if $N>2$. However, if $N=2$, then $\bar {\alpha}=2/p'$, $\gamma_p(\bar{\alpha},0)=0$,  hence Rellich inequalitiy holds for $\bar {\alpha}$ if and only if $N\ge 3$, according to (\ref{condRp}).
Assume now that $\alpha \neq \bar{\alpha}$.
Since $\lambda_n=n^2+(N-2)n$, it follows from (\ref{gammap}) with $c=0$ that
$$
\lambda_n+\gamma_p(\alpha,0)=\Bigl(n+\frac{N}{2}-1\Bigr)^2-\Bigl(N\Bigl (\frac12-\frac1p\Bigr)+1-\alpha\Bigr)^2=\left (\frac{N}{p}-2+\alpha+n\right )\left (\frac{N}{p'}-\alpha+n\right ),
$$
hence the condition $\lambda_n+\gamma_p(\alpha,0) \neq 0$ for every $n \in \N_0$ translates into (\ref{condRp}). Finally, if $2-N/p<\alpha <N/p'$, then $\gamma_p(\alpha,0)>0$ and the best constant is given by (\ref{bestp}), by Theorem \ref{Rellich-p}.
\qed

When Rellich inequalities hold for the Laplacian other inequalities of Sobolev type can be proved. We refer the reader to the very recent paper \cite{musina} where this topic is studied sistematically and confine ourselves to quote the following result, see \cite[Corollary 2.12 (i)]{musina}. If $p <N/2$,  $p^{**}$ is defined by $1/p^{**}=1/p-2/N$.

\begin{prop} \label{wsobolev}  Let $N \ge 2$, $1< p< \infty$, $\alpha \in \R$, $b=c=0$ and assume that (\ref{condRp}) holds. Then for $p \le q \le p^{**}$ when $p<N/2$ and and for  $p \le q <\infty$ for $p \ge N/2$ there exists $C>0$ such that
$$
C\int_{\R^N}|x|^{-N+q\frac{N-2p+\alpha p}{p}}|u|^q\, dx \le \int_{\R^N}|x|^{\alpha p}|\Delta u|^p\, dx
$$
for every $u\in C_c^\infty(\R^N \setminus \{0\})$. 
\end{prop}

The two conditions in (\ref{condRp}) are not independent and the best constants $C(N,\alpha,p,0,0)$ satisfy a simmetry relation.
\begin{cor}
Let $N \ge 3$, $\alpha \in \R$ and  $\beta=2-\alpha +N(1-2/p)$. Rellich inequalities hold for the Laplacian (i.e. with b=c=0) for the weight $|x|^\alpha$ if and only if they hold for the weight $|x|^\beta$. Moreover the best constants satisfy $C(N,\alpha,p,0,0)=C(N,\beta,p,0,0)$.
\end{cor}
{\sc Proof. }  We use the Kelvin transform $u(x)=|x|^{2-N}v\left (\frac{x}{|x|^2}\right )$ when $N \ge 3$ and $u,v \in C_c^\infty (\R^N \setminus \{0\})$. Then 
$$\Delta u(x)=|x|^{-N-2}\Delta v\left (\frac{x}{|x|^2} \right ).$$
Setting $y=x/|x|^2$, $dx=|y|^{-2N} dy$ and  by elementary computations we see that the inequality
$$
\||x|^\alpha \Delta u\|_p \geq C\||x|^{\alpha-2} |u|\|_p
$$
is equivalent to 
$$
\||x|^\beta \Delta u\|_p \geq C\||x|^{\beta-2} |u|\|_p
$$
with the same constant $C$. \qed

Note that $\alpha \neq N/p'+n$ is equivalent to $\beta \neq -N/p+2-n$

\begin{os}
{\rm The computation of the best constant when $2-N/p<\alpha <N/p'$ (which requires $N \ge 3$) is due to Mitidieri, \cite[Theorem 3.1]{mitidieri}, and also to Davies and Hinz, \cite[Theorem 12]{davi-hinz},  under the more restrictive condition $2-\frac{N}{p}<\alpha<2-\frac{2}{p}$. When  $\alpha=0$ and $1<p<N/2$ the best constant is given $
\gamma_p(0,0)=\Bigl(\frac{N}{p}-2\Bigr)\frac{N}{p'}>0$, according to \cite{okazawa}, as mentioned in the Introduction.
}
\end{os}

\begin{os}
{\rm It will be shown in Section 7 (see Theorem \ref{RellichFJ}) that when one of the conditions in (\ref{condRp}) is violated for a specific $n \in \N_0$, then Rellich inequalities fail for functions of the form $\sum_jf_j(\rho)P_j(\omega)$, where $f_j $ are  smooth functions with compact support and  $P_j$ are spherical harmonics of order $n$. However Rellich inequalities hold in the "complementary subspace", that is for functions $u=\sum_k g_k(\rho)P_k(\omega)$, where $g_k$ are as above and $P_k$ are spherical harmonics of order different from $n$.
}
\end{os}

\begin{os}
{\rm By inverting the role of $p, \alpha$ in (\ref{condRp}) we may identify, for a fixed $\alpha$, the values of $p$ for which Rellich inequalities fail. It follows that
\begin{itemize}
\item[(i)] If $\alpha \ge 0$ Rellich inequalities fail for all $p=N(N+n-\alpha)$ if $(\alpha-N)\vee 0 \le n \le \alpha$.
\item[(ii)] If $\alpha \le 2$ Rellich inequalities fail for all $p=N/(2-n-\alpha)$ if $((2-\alpha)-N)\vee 0 \le n \le 2-\alpha$.
\end{itemize}
Observe that there is an overlapping between (i) and (ii) on the interval $[0,2]$. In particular, if $\alpha=0$, Rellich inequalities fail for $p=1,N/2, N, \infty$. If $p=1,N/2$ they fail for $n=0$, that is for radial functions,  but hold for smooth functions $u$ having zero mean on $S_{N-1}$, that is 
$$\int_{S_{N-1}}u(\rho, \omega) \, d\sigma(\omega)=0.$$ To see this it is sufficient to apply  the next Theorem \ref{RellichFJ} with $J=\{0\}$ and Lemma \ref{projection}.
Rellich inequalities fail  for $ p=N$ and $n=1$, that is when the spherical harmonics have order 1, but  hold for smooth functions $u $ such that
$$\int_{S_{N-1}}u(\rho, \omega)P(\omega) \, d\sigma(\omega)=0$$
for every spherical harmonic of order 1 (apply   Theorem \ref{Cald-ZygFJ} with $J$ corresponding to the spherical harmonics of order 1 and Lemma \ref{projection}. Finally, Rellich inequalities fail for $p=\infty$ with $n=2$.
However Rellich inequalities hold for every $1\le p \le \infty$ for functions $u=\sum_k g_k(\rho)P_k(\omega)$, if the spherical harmonics $P_k$ have order at  least 3. We refer again to Section 7 and, in particular, to Theorem \ref{RellichFJ}.
}
\end{os}

\subsection{Complex coefficients}
First we observe that the proof of Theorem \ref{Rellich-p} works also when $b$ is complex, except for the computation of the optimal constant.
The case where also $c$ is complex can be easily reduced to the previous one.
Let

\begin{equation} \label{gammapC}
\gamma_p(\alpha,c)=\Bigl(\frac{N}{p}-2+\alpha+ i\frac{c_2}{2}\Bigr)\Bigl(\frac{N}{p'}-\alpha+c_1+i\frac{c_2}{2}\Bigr)=\Bigl(\frac{N}{2}-1+\frac{c}{2}\Bigr)^2-\Bigl(N\Bigl (\frac12-\frac1p\Bigr)+1-\alpha+\frac{c_1}{2}\Bigr)^2
\end{equation}
and the parabola
\begin{equation} \label{PpC}
{\cal P}_{p,\alpha,c}:=
\left\{\lambda=-\xi^2+i\xi\Bigl (N\Bigl (1-\frac2p \Bigr)+2-2\alpha+c_1\Bigr)-\gamma_p(\alpha,c)\;;\;\xi\in \R\right\}.
\end{equation}

\begin{prop} 
Let $N\geq 2$, $\alpha \in \R$, $b,\ c=c_1+ic_2\in\C$, $1\le p\le \infty$.  There exists a positive constant $C=C(N,\alpha,  p, c, b)$ such that
\begin{equation} \label{rellich-op-comp}
\||x|^\alpha Lu\|_p \geq C\||x|^{\alpha-2} |u|\|_p
\end{equation}
for every $u\in C_c^\infty(\R^N \setminus \{0\})$ if and only if  $b+\lambda_n \not \in {\cal P}_{p,\alpha,c}$ for every $n \in \N_0$.
 If, in addition, $\rm Im\, b+\rm Im\,\gamma_p(\alpha,c)=0$ and $\rm Re\, b+\rm Re\, \gamma_p(\alpha,c)>0 $
the optimal constant is given by
$C=\rm Re\, b+\rm Re\, \gamma_p(\alpha,c).$
 \end{prop}
{\sc Proof.} We have already observed that the proof of Theorem \ref{Rellich-p} holds if $b \in \C$ and $c \in \R$. Consider now the general case $c,\ b\in\C$ with $c=c_1+ic_2$.
Let $u\in C_c^\infty(\R^N\setminus\{0\})$. As before, setting $v(x)=|x|^{\alpha-2}u(x)$, we have
$$|x|^\alpha L u:=\tilde{L}v-bv$$
where 
\begin{equation*}
\tilde L=|x|^2\Delta +(4-2\alpha+c)x\cdot\nabla +(2-\alpha)(N-\alpha+c).
\end{equation*}
A second transformation allows us to obtain an operator with real drift.
Set $$T: C_c^\infty(\R^N\setminus\{0\})\to C_c^\infty(\R^N\setminus\{0\}),\quad u\to |x|^{i\gamma}u,$$ with $\gamma=-\frac{c_2}{2}$.
Observe that $T$ is an isometry in $L^p(\R^N)$ and 
$$T^{-1}\tilde{L}T=|x|^2\Delta +(4-2\alpha+c_1)x\cdot\nabla-\left[i\frac{c_2}{2}\left(i\frac{c_2}{2}+N-2+c_1\right)-(2-\alpha)(N-\alpha+c_1)\right].$$
Therefore estimate (\ref{rellich-op-comp}) is equivalent to the estimate
\begin{equation*} 
\|T^{-1}(\tilde{L}-b)Tv\|_p=\|T^{-1}\tilde{L}Tv-bv\|_p \geq C\|v\|_p
\end{equation*}
for any $v\in C_c^\infty (\R^N \setminus \{0\})$.  The last estimate is true if and only if $b$ does not belong to the  spectrum of  the operator 
$T^{-1}\tilde{L}T$ which  is, by Proposition \ref{specRN},
$$\sigma(T^{-1}\tilde{L}T)=\cup_{n\in\N_0}({\cal P}_{p,\alpha,c}-\lambda_n).$$
The optimal constant can be computed as before in the case $\rm Im\, b+\rm Im\,\gamma_p(\alpha,c)=0$ and $\rm Re\, b+\rm Re\, \gamma_p(\alpha,c)>0$.
\qed

\subsection{One dimension}
The results in the previous sections have been stated and proved for $N \ge 2$. However they also holds in one dimension with similar but simpler proofs. We formulate the next result in $]0,\infty[$; the case of the whole space follows immediately by adding the corresponding inequalities in $]-\infty,0[$ and $]0,\infty[$.
According with the previous notation, we set 
$$L=D^2 +c\frac{x}{|x|^2}D-\frac{b}{|x|^2},$$ 
\begin{equation*} 
\gamma_p(\alpha,c)=\Bigl(\frac{1}{p}-2+\alpha\Bigr)\Bigl(\frac{1}{p'}-\alpha+c\Bigr)
\end{equation*}
and 
\begin{equation} \label{Pp1}
{\cal P}_{p,\alpha,c}:=
\left\{\lambda=-\xi^2+i\xi\Bigl (2-2\alpha+c\Bigr)-\gamma_p(\alpha,c)\;;\;\xi\in \R\right\}.
\end{equation}
For simplicity we assume that $c \in \R$.

\begin{prop} \label{Rellich-p1}
Let $N=1$, $\alpha,\ c\in\R$, $b \in \C$,  $1 \le p \le \infty$.
Then there exists a positive constant $C=C(\alpha,  p, c, b)$ such that
\begin{equation} \label{rellich-op1}
\||x|^\alpha Lu\|_p \geq C\||x|^{\alpha-2} |u|\|_p
\end{equation}
holds for every $u\in C_c^\infty(]0,\infty[)$ if and only if $b \not \in {\cal P}_{p,\alpha,c}$.
If, in addition, $b \in \R$ and  $b+\gamma_p(\alpha,c)>0$
the optimal constant is given by
$C=b+\gamma_p(\alpha,c).$
\end{prop}

The statement follows arguing as in the case $N\geq 2$. In this case the auxiliary operator $\tilde L$ is one-dimensional and neither spherical harmonics nor eigenvalues $\lambda_n$ appear.

\section{Weighted Calder\'on-Zygmund inequalities}

Rellich inequalities can be used to compute all $\alpha's$ for which weighted Calder\'on-Zygmund inequalities hold in the weaker form (\ref{cald-zygmund}) below. We refer to  the Introduction for a comparison between 
(\ref{cald-zygmund}) and the stronger form (\ref{Stein}).

\begin{teo} \label{Cald-Zyg}
Let $N\geq 3$, $1 < p < \infty$, $\alpha\in\R$.
The
weighted Calder\'on-Zygmund inequalities
 \begin{equation} \label{cald-zygmund}
\||x|^\alpha D^2 u\|_p \leq C \||x|^{\alpha} \Delta u\|_p
\end{equation}
hold in $C_c^\infty(\R^N \setminus \{0\})$ if and only if  

$$
\alpha \neq \frac{N}{p'}+n \quad {\rm for\  every\ } n \ge 0 \quad {\rm and\  }\quad \alpha \neq -\frac{N}{p}+2-n \quad {\rm for\  every\ } n \ge2.
$$
\end{teo}

\begin{os}
{\rm If $N=2$, $\ds\frac{N}{p'}$ coincides with $2-\frac{N}{p}$. In such a case both weighted Calder\'on-Zygmund and Rellich inequalities fail. 
Indeed, the family of functions $\{u_m\}_{m\in \N}$ defined as
\[
u_m(x)= m^{2-\frac{1}{p}}\phi\left(\frac{\log|x|}{m}\right),
\]
where $\phi\in C_c^\infty(]0,+\infty[)\setminus\{0\}$, 
satisfies that $\||x|^{\alpha}\Delta u_m\|_p$ is independent of $m$ 
and $\||x|^{\alpha}D^2u_m\|_p\to \infty$ as $m\to \infty$. 
More precisely, we have 
\begin{align*}
\||x|^{\alpha}\Delta u_m\|_p=(2\pi)^{\frac{1}{p}}\|\phi''\|_{L^p(]0,+\infty[)}.
\end{align*}
On the other hand, we see that 
\begin{align*}
\||x|^{\alpha}D^2 u_m\|_p
&\geq 
\frac{1}{4}\left(
    \frac{\pi}{6}
\right)^{\frac{1}{p}}( 2m \|\phi'\|_{L^p(]0,+\infty[)}-\|\phi''\|_{L^p(]0,+\infty[)}). 
\end{align*}
Hence 
we conclude that the sequence $\{u_m\}_m$ is a counterexample of weighted Calder\'on-Zygmund inequalities.}
\end{os}

We need some preliminary interpolative estimates.

\begin{lem} \label{inter1d}
Let  $1 \le p \le \infty$ and $\beta \in \R$. Then there exists $C=C(p,\beta)>0$ such that for every $u\in C_c^\infty (]0,\infty[)$
$$
\|\rho^{\beta-1} u'\|_p \le \eps \|\rho^\beta u''\|+\frac{C}{\eps}\|\rho^{\beta-2}u\|_p,
$$
for $0<\eps \le 1$, the norms being taken on $(0,\infty)$.
\end{lem}
{\sc Proof.} By Taylor's formula and for $h>0$
$$
u(\rho+h)-u(x)=hu'(\rho)+\int_0^h (h-s)u''(\rho+s)\, ds.
$$
Setting $h=\eps \rho$ we obtain
$$
u'(\rho)=\frac{u((1+\eps)\rho)-u(\rho)}{\eps \rho}+\eps \int_0^1 \rho(1-h)u''(\rho(1+\eps h))\, dh
$$
hence
$$\rho^{\beta-1}u'(\rho)=\eps^{-1}\rho^{\beta-2}\left (u((1+\eps)\rho)-u(\rho)\right )+\eps \int_0^1 \rho^\beta(1-h)u''(\rho(1+\eps h))\, dh
$$
Taking the $L^p$-norms of both sides and using Minkowski inequality for integrals, the result follows by easy computations (note that all integrals with respect to the $\rho$ variable are uniformly bounded in $0<\eps,h \le 1$).
\qed

Nex we prove the N-dimensional version of the above lemma.

\begin{lem}  \label{interpolative}
Let $\alpha\in\R$, $1 \le p \le \infty$. There exist $C(N, p,\alpha)>0$, $\eps_0>0$  such that for every $0<\varepsilon\leq \eps_0$ and any $u\in C_c^\infty(\R^N\setminus\{0\})$,
\begin{equation}\label{eq:stima-interpolativa}
\||x|^{\alpha -1}\nabla u\|_p
\leq \varepsilon \||x|^\alpha D^2 u\|_p
    +\frac{C}{\varepsilon} \||x|^{\alpha-2} u \|_p.
\end{equation}
\end{lem}
{\sc Proof.} Let $u=u(\rho\, \omega)$ with $|\omega|=1$. Then $u_\rho=\nabla u \cdot \omega$, $u_{\rho \rho}=\sum_{i,j}D_{ij}u\  \omega_i \omega_j$. We apply Lemma \ref{inter1d} with $\beta=\alpha+(N-1)/p$ and obtain
$$
\int_0^\infty \rho^{p(\alpha-1)+N-1}|u_\rho|^p\, d\rho \le \eps^p \int_0^\infty \rho^{p\alpha+N-1}|u_{\rho \rho}|^p\, d\rho +\frac{C^p}{\eps^p}\int_0^\infty \rho^{p(\alpha-2)+N-1}|u|^p\, d\rho.
$$ 
Integrating the above inequality with respect to $\omega \in S_{N-1}$ we obtain 
\begin {equation} \label{radiale}
\int_{R^N} |x|^{p(\alpha-1)}|u_\rho|^p\, dx \le \eps^p \int_{\R^N}|x|^{p(\alpha-2)}|D^2 u|^p\, dx+\frac{C^p}{\eps^p }\int_{R^N}|x|^{p(\alpha-2)}|u|^p\, dx,
\end{equation}
that is (\ref{eq:stima-interpolativa}) for the radial component of the gradient.
Concerning the analogous estimate for the tangential gradient we observe that if $v \in C^\infty (S_{N-1})$ then the classical interpolative estimate
$$
\int_{S_{N-1}}|\nabla_{\tau} v|^p\, d\sigma \le \eps^p \int_{S_{N-1}}|D^2_{\tau} v|^p\, d\sigma +\frac{C^p}{\eps^p}\int_{S_{N-1}}|v|^p\, d\sigma 
$$
holds, where $\nabla_{\tau}$ and $D^2_{\tau}$ denote the tangential gradient and the tangential Hessian matrix, respectively.
Applying it to $v(\omega)=u(\rho\,  \omega)$, multiplying by $\rho^{p(\alpha-2)+N-1}$ and integrating over $(0, \infty)$ we obatin we obtain
\begin {equation} \label{tangenziale}
\int_{R^N} |x|^{p(\alpha-2)}|\nabla_{\tau}u|^p\, dx \le \eps^p \int_{\R^N}|x|^{p(\alpha-2)}|D^2_{\tau} u|^p\, dx+\frac{C^p}{\eps^p }\int_{R^N}|x|^{p(\alpha-2)}|u|^p\, dx.
\end{equation}
Since $\nabla u=u_\rho \frac{x}{\rho}+\frac{1}{\rho}\nabla_\tau u$, $|\nabla u|^2=u_\rho^2+\frac{1}{\rho^2}|\nabla_\tau u|^2$ and since $|D^2_\tau u|$ is pointwise dominated by $\rho^2|D^2u|+\rho |\nabla u|+|u|$, suming (\ref{radiale}) and (\ref{tangenziale}) and taking $\eps$ small we conclude the proof.
\qed

\begin{os}
{\rm Clearly Lemma \ref{interpolative} holds for functions $u \in W^{2,p}$ having compact support in $\R^N \setminus \{0\}$. Moreover it holds in subspaces of $Z \subset W^{2,p}(\Omega)$ for which there exists a linear extension operator $E$ from $Z$ to $W^{2,p}$ functions with compact support in $\R^N\setminus \{0\}$.}
\end{os}

{\sc Proof of Theorem \ref{Cald-Zyg}}
We first show that Rellich inequalities imply Calder\'{o}n-Zygmund inequalities. We  apply the classical
Calder\'{o}n-Zygmund inequality
$$
\|D^2u \|_p \le C\|\Delta u\|_p
$$
to $|x|^\alpha u$ and 
estimate the first order terms using  Lemma \ref{interpolative}. We get for small $\eps$
\begin{align*}
\||x|^\alpha D^2 u\|_p&\leq
  C\left( \|D^2(|x|^\alpha u)\|_p+\||x|^{\alpha-1}\nabla u\|_p +\||x|^{\alpha-2}u\|_p\right)\\
  &\leq
C\left( \|\Delta(|x|^\alpha
u)\|_p+\||x|^{\alpha-1}\nabla u\|_p +\||x|^{\alpha-2}u\|_p\right)\\&\leq
C\left( \||x|^\alpha\Delta
u\|_p+\eps \||x|^{\alpha}D^2 u\|_p +C_\eps\||x|^{\alpha-2}u\|_p\right).
\end{align*}
Taking $\eps$ such that $C\eps <1/2$ and by applying Rellich inequalities, weighted Calder\'{o}n-Zygmund inequalities follow.

By Theorem \ref{CMp} we obtain that 
Calder\'{o}n-Zygmund inequalities hold when $\alpha \neq N/p' +n$, $\alpha \neq -N/p+2-n$, $n \in \N_0$.
However,  if $n=0,1$ in the second formula, that is if  $\alpha=2-N/p$, $\alpha=1-N/p$,  then Calder\'{o}n-Zygmund inequalities hold by Stein result \cite{stein}, since $-N/p <\alpha<N/p'$ (here we need $N \ge 3$).

\smallskip

\noindent Let us now assume that Calder\'{o}n- Zygmund inequalities hold and that 
$$\alpha \neq 2-\frac{N}{p},\quad \alpha \neq 1-\frac{N}{p}.$$

We may therefore  apply Hardy inequalities twice (use Proposition \ref{hardy} with $\beta=(\alpha-1)p$ and $\beta=\alpha p$) to obtain 
 \begin{equation} \label{hardytwice}
\left |\left (\frac{N}{p}+\alpha-2\right )\left (\frac{N}{p}+\alpha-1 \right )\right | \||x|^{\alpha-2} u\|_p \leq  \||x|^{\alpha} D^2 u\|_p
\end{equation}
 for every $u\in C_c^\infty(\R^N \setminus \{0\})$.
By the weighted Calder\'on-Zygmund inequalities, Rellich inequalities follows and then  
$$
\alpha \neq \frac{N}{p'}+n \quad \alpha \neq -\frac{N}{p}+2-n
$$
for $n \in \N_0$,
by Theorem \ref{CMp}.
 \qed

\begin{os}
{\rm Observe that Rellich inequalities do not hold  when  $\alpha= 1-\frac{N}{p}$,  $\alpha= 2-\frac{N}{p}$ but  Calder\'{o}n-Zygmund inequalities are true.}
\end{os}

\begin{os}
{\rm  Calder\'{o}n-Zygmund inequalities have been used in \cite{met-spi1}, \cite{met-spi2}, \cite{met-spi-tac} to characterize the domain of second order elliptic operator with unbounded coefficients like $|x|^\alpha \Delta$ or $(1+|x|^\alpha )\Delta +c|x|^{\alpha-1}\frac{x}{|x|}\cdot \nabla$. Even though (\ref{cald-zygmund}) holds in most cases, the characterization of the domain is possible under more restrictive conditions ensuring the density of smooth functions (where (\ref{cald-zygmund}) holds) in the domain of the operator. In order to describe the domain, the weaker inequality
\begin{equation} \label{CZw}
\||x|^\alpha D^2 u\|_p \leq C (\||x|^{\alpha} \Delta u\|_p+\|u\|_p)
\end{equation} 
 suffices. However, if $\alpha \neq 2$ this weaker inequality implies the stronger (\ref{cald-zygmund}) by replacing $x$ with $\lambda x$ and then letting $\lambda \to 0, \infty$ (according to the sign of $\alpha-2$). This argument fails if $\alpha=2$ and, in fact, (\ref{CZw}) always holds for $\alpha=2$ and every $1<p<\infty$, see \cite{met-spi1}, but (\ref{cald-zygmund}) fails for $N \ge 3$ and $p=N/(N-2)$.
}
\end{os}

Theorem \ref{Cald-Zyg} can be partially generalized to the case of more general operators $L=\Delta +c\frac{x}{|x|^2}\cdot \nabla -\frac{b}{|x|^2}$.

\begin{prop} \label{CZL}
Let $1<p<\infty$, $\alpha \in \R$ and assume that Rellich inequalities hold for $L$. Then 
Calder\'{o}n-Zygmund inequalities
$$
\||x|^\alpha D^2 u\|_p \le C\||x|^\alpha Lu\|_p
$$
hold in $C_c^\infty (\R^N\setminus \{0\})$ with a suitable $C>0$.
\end{prop}
{\sc Proof. } As in the proof of Theorem \ref{Cald-Zyg} we obtain
\begin{align*}
\||x|^\alpha D^2 u\|_p&\leq
  C\left( \||x|^\alpha \Delta u\|_p+\||x|^{\alpha-1}\nabla u\|_p +\||x|^{\alpha-2}u\|_p\right)\\
  &\leq
C\left( \|x|^\alpha L
u\|_p+\||x|^{\alpha-1}\nabla u\|_p +\||x|^{\alpha-2}u\|_p\right)\\&\leq
C\left( \||x|^\alpha L
u\|_p+\eps \||x|^{\alpha}D^2 u\|_p +C_\eps\||x|^{\alpha-2}u\|_p\right).
\end{align*}
Taking $\eps$ such that $C\eps <1/2$ and by applying Rellich inequalities, weighted Calder\'{o}n-Zygmund inequalities follow.
\qed
Calder\'{o}n-Zygmund inequalities hold for $L$ whenever $b+\lambda_n \not \in {\cal P}_{p,\alpha,c}$ for every $n \in \N_0$, in particular for $b=0$, $2 \le \alpha <N/p'+c $, a case first established in \cite{met-spi-tac} in the framework of elliptic operators with unbounded drift and diffusion coefficients. 

\section{Spectrum of the operator $A=|x|^2\Delta+ cx\cdot\nabla$} \label{spectrumA}
In this section we compute the spectrum of the operator 
$$A=|x|^2\Delta+ cx\cdot\nabla$$ in $L^p({\cal C}_\Sigma)$, where $\Sigma \subseteq S_{N-1}$ is relatively open and $C^2$ and 
$${\cal C}_\Sigma=\{x=(\rho,\omega)\in\R^N: \quad \rho>0,\ \omega\in\Sigma\}$$
is the cone with vertex at $0$ defined by $\Sigma$. We are mainly interested in the cases of the whole space, corresponding to $\Sigma=S_{N-1}$ and of the half-space, corresponding to $\Sigma=S_{N-1}^+=\{\omega \in S_{N-1}: \omega_N>0\}$. When $\Sigma \not =S_{N-1}$ we impose Dirichlet boundary conditions on  $\partial C_\Sigma \setminus \{0\}$.
In order to compute the spectrum we write $A$ in spherical coordinates 

$$A=\rho^2\frac{\partial^2}{\partial\rho^2}+(N-1+c)\rho \frac{\partial }{\partial \rho}+\Delta_0,$$

\noindent
 and consider the operators 
\begin{equation} \label{opgamma}
\Gamma=\rho^2\frac{\partial^2}{\partial\rho^2}+(N-1+c)\rho \frac {\partial }{\partial \rho}
\end{equation}
 in $L^p((0,\infty),\ r^{N-1}dr)$ and in $C_0^0(]0,+\infty[)$ (consisting of all continuous functions vanishing at $0,\infty$) and the Laplace-Beltrami operator 
$\Delta_0$ in $L^p(\Sigma)$ and in $C_0(\Sigma)$, endowed with Dirichlet boundary conditions on $\partial \Sigma$, the boundary of $\Sigma$ in $S_{N-1}$.
The operators $\Gamma$ and $\Delta_0$ act on independent variables and therefore the spectrum of their sum can be computed through tensor products arguments. We start by analyzing them separately.
\subsection{The operator $\Gamma$}

To shorten the notation, when $p=\infty$, $L^p(]0,+\infty[,\rho^{N-1}\,d\rho)$ stands for $C_0^0(]0,+\infty[)$.
The operator $\Gamma$ is defined in (\ref{opgamma}).

\begin{prop} \label{spec-rad}
Let $1\leq p\leq \infty$. Then $\Gamma$, endowed with the domain $$D_p(\Gamma)=\{u\in L^p(]0,\infty[,\rho^{N-1}\,d\rho),\ \rho\frac{\partial u}{\partial \rho},\ \rho^2\frac{\partial^2 u}{\partial \rho^2}\in L^p(]0,\infty[,\rho^{N-1}\,d\rho)\},$$  generates a strongly continuous and analytic semigroup $(S(t))_{t \ge 0}$  in $L^p(]0,\infty[, \ \rho^{N-1}d\rho)$.\\ Its spectrum is given by 
\begin{equation} \label{spettrogamma}
\sigma_p(\Gamma)={\cal P}_p=\left\{\lambda=-\xi^2+i\xi\left (N(1-\frac{2}{p})-2+c\right)-\omega_p,\,  \xi\in\R\right\},
\end{equation}
where  
\begin{equation} \label {defomegap}
\omega_p=\frac{N}{p^2}\left[p(N-2+c)-N\right]
\end{equation}
and  $(S(t))_{t \ge 0}$ satisfies the estimate $\|S(t)\|_p \le e^{-\omega_p t}$ for $t \ge 0$.
\end{prop}
{\sc Proof.} 
Consider the transformations $$S:L^p(\R, ds)\to L^p(]0,\infty[),\rho^{N-1}\,d\rho), \quad (Su)(r)=\rho^{-\frac{N}{p}}u(\log \rho),$$ for $1\leq p<\infty$, and   $$S:C_0(\R)\to C_0^0(]0,\infty[), \quad Su(r)=u(log \rho).
$$ It is easy to show that $S$ is an isometry and that
$$S^{-1}\Gamma S u=u''+\left(N-2-\frac{2N}{p}+c\right)u'-\omega_p u$$ hence,
by classical results,  $(S^{-1}\Gamma S, W^{2,p}(\R) )$  and  $(S^{-1}\Gamma S, C_0^2(\R) )$ generate a strongly continuous analytic semigroup in $L^p(\R)$  and  in $C_0(\R)$, respectively, whose norm is bounded by $e^{-\omega_p t}$. It follows that $\Gamma$, endowed with the domains
$$D_p(\Gamma)=\{Su:\ u\in W^{2,p}(\R)\} \subset L^p(]0,\infty[, \rho^{N-1}d\rho[)$$
and
$$D_\infty(\Gamma)=\{u(\log \rho):\ u\in C_0^2(\R)\} \subset C_0^0(]0,\infty[)$$ generates a strongly continuous and analytic semigroup $(S(t))_{t \ge 0}$  in $L^p(]0,\infty[, \ \rho^{N-1}d\rho)$ satisfying   $\|S(t)\|_p \le e^{-\omega_p t}$. 
It is easy to check that
$$D_p(\Gamma)=\{u\in L^p(]0,\infty[,\rho^{N-1}\,d\rho),\ \rho\frac{\partial u}{\partial \rho},\ \rho^2\frac{\partial^2 u}{\partial \rho^2}\in L^p(]0,\infty[,\rho^{N-1}\,d\rho)\},$$for $1\leq p<\infty$, and  
$$D_\infty(\Gamma)=\{u\in C_0^0(]0,\infty[):\ \rho\frac{\partial u}{\partial \rho},\ \rho^2\frac{\partial^2 u}{\partial \rho^2}\in C_0(]0,\infty[)\}$$
in $C_0^0(0,+\infty)$. \\

\noindent
Concerning the second part of the statement we observe that the spectra of $\Gamma$ and $S^{-1}\Gamma S$ coincide. The operator $S^{-1}\Gamma S$  is uniformly elliptic in $L^p(\R, ds)$, hence  its spectrum is independent of $p$ and coincides with the spectrum in $L^2(\R, ds)$ which can be computed using the Fourier transform ${\cal F}$.
Since $${\cal F}{(S^{-1}\Gamma S)}u(\xi)=\left(-\xi^2+i\xi\left(N-2+c-\frac{2N}{p}\right)-\omega_p\right){\cal F}{u},$$ the formula for ${\cal P}_p$ follows.
\qed

In order to compute best constants in some Rellich inequalities we need the norm of the resolvent of $\Gamma$ for $\lambda \in \R, \lambda \not \in {\cal P}_p$. 

\begin{lem} \label{norma}
Let $B=D^2+2b D$ in $L^p(\R)$, $1\le p\le \infty$. If $\lambda \in \R$, $\lambda \not =0$ and $\lambda \ge -b^2$, then 
$\|(\lambda-B)^{-1}\|_p=|\lambda|^{-1}$.
\end{lem}
{\sc Proof. } The spectrum of $B$ is independent of $p$ and given by the parabola
$$
{\cal P}=\{\lambda=-\xi^2+2ib\xi, \, \xi \in \R\}
$$
and therefore $\|(\lambda-B)^{-1}\|_p \ge dist^{-1} (\lambda, {\cal P})\ge  |\lambda|^{-1}$ for $\lambda \not \in {\cal P}$.
Next observe that $\|(\lambda-B)^{-1}\|_p \le \lambda^{-1}$ for $\lambda >0$, since $B$ generates a contraction semigroup in $L^p(\R)$. The assertion is then proved for $\lambda>0$. \\

\noindent
Assume now that $-b^2 \le \lambda <0$. 
We write explicitly the resolvent and assume for example that $b>0$. For $-b^2 <\lambda<0$ two linearly independent solutions of the homogenuous equation $\lambda u-Bu=0$ are given by $u_i(t)=e^{\mu_i t}$, $ \mu_i=-b\pm \sqrt{b^2+\lambda}$,  $i=1,2$. Variation of constants yields $u=T_Kf$ if $\lambda u-Bu=f$, where
$$
T_Kf(t)=\int_{-\infty}^{+\infty} K(t,s)f(s)\, ds
$$
and
$$
K(t,s)=\frac{1}{\mu_2-\mu_1}\left (e^{\mu_1(t-s)}-e^{\mu_2 (t-s)} \right )\chi_{\{s\le t\}}
$$
(observe that $K$ is negative). Therefore
$$
\|(\lambda-B)^{-1}\|_p \le \max \left\{\sup_t \int_{-\infty}^{+\infty}|K(t,s)|\, ds, \sup_s \int_{-\infty}^{+\infty}|K(t,s)|\, dt \right \}=|\lambda|^{-1}.
$$
This gives the result for $-b^2 <\lambda <0$ and, by continuity, also for $\lambda=-b^2$, when $b \not =0$.
\qed

\begin{prop} \label{normagamma} Assume that $ \lambda \in \R$ and that $2(\lambda-\omega_p)  +\left (N(1-\frac{2}{p})-2+c\right)^2\ge 0$. Then 
$\|(\lambda-\Gamma)^{-1}\|_p=|\lambda-\omega_p|^{-1}$. 
\end{prop}
{\sc Proof.} As in Proposition \ref{spec-rad}
$$S^{-1}\Gamma S =D^2+\left(N-2-\frac{2N}{p}+c\right)D-\omega_p $$ and $S$ is an isometry of the corresponding spaces. Therefore the thesis follows from Lemma \ref{norma}.
\qed

\begin{os}
{\rm The equality $\|(\lambda-B)^{-1}\|_p=|\lambda|^{-1}$ stated in Lemma \ref{norma} is true in $L^2(\R)$ if and only if $\lambda \ge -2b^2$, since it coincides with the distance of $\lambda$ from the parabola ${\cal P}$, see Lemma \ref{dist}. The resolvent of $B$ can be computed also for $\lambda<-b^2$. Writing $\lambda=-b^2-\gamma^2$ we obtain
$(\lambda-B)^{-1}f=f*g$ where $g(t)=\gamma^{-1}e^{-b t}\sin (\gamma t)\chi_{\{t \ge 0\}}$. Hence, putting $s=\gamma t$ and summing the integrals where $\sin s$ is positive and negative,
$$
\|(\lambda-B)^{-1}\|=\frac{1}{\gamma}\int_0^\infty e^{-b t}|\sin (\gamma t)|\, dt=\frac{1}{b^2+\gamma^2}\coth \frac{b\pi}{2\gamma}
$$
for $p=1, \infty$. The norm for other values of $p$ can be estimated by interpolating between $p=1,2,\infty$.
}
\end{os}

\subsection{The operator $\Delta_0$}
The Laplace Beltrami operator $\Delta_0$, endowed with Dirichelet boundary conditions (if $\Sigma \not = S_{N-1}$), generates an analytic semigroup $(T_\Sigma(t))_{t \ge 0}$ in $L^p(\Sigma)$ (with respect to the surface measure $d\sigma$) for every $1\leq  p < \infty$ and in $C_0(\Sigma)$. By elliptic regularity it follows that its domain $D_p(\Delta_0, \Sigma)$ is given by $W^{2,p}(\Sigma, d\sigma)\cap W^{1,p}_0(\Sigma, d\sigma)$ if $1<p<\infty$. The analyticity of the semigroup follows from  Gaussian estimates of the  heat kernel of $\Delta_0$ proved in \cite[Theorem 5.2.1, Theorem 5.5.1]{davies}, using \cite[Corollary 7.5]{ou}.
\smallskip\noindent

\begin{lem} \label{eigen}
The spectrum of the operator $(\Delta_0,D_p(\Delta_0, \Sigma))$ is independent of $1 \le p \le \infty$ and consists of isolated eigenvalues. Each eigenvalue is a simple pole of the resolvent and has a finite geometric multiplicity which is equal to its algebraic multiplicity. The eigenfunctions are  independent of $p$ and their linear span is dense in $L^p(\Sigma)$ for $1 \le p<\infty$ and $w^*$-dense for $p=\infty$.
\end{lem}
{\sc Proof.} The operator is self-adjoint in $L^2(\Sigma)$ and has a compact resolvent for every $1 \le p \le \infty$. The independence of the spectrum, as well as of the spectral projections, multiplicities and of the eigenfunctions follows from classical results, see e.g. \cite[Proposition 2.2]{arendt}. Each eigenvalue is a simple pole of the resolvent for $p=2$ since the operator is self-adjoint and hence for every $p$, since the Laurent expansion of the resolvent around each eigenvalue is independent of $p$. The equality of geometric and algebraic multiplicities follows from self-adjointness in $L^2$, and then in the general case since these quantities are independent of $p$ (see also \cite[Proposition 5.5]{FMPP} where these arguments are explained in more detail).
Let $P_n(\Sigma)$ be the eigenfunctions and note that they belong to $L^p(\Sigma)$ for every $p$ and that they form a complete orthonormal system in $L^2(\Sigma)$ since the operator is self-adjoint. First let $p<2$ and consider $f \in L^2(\Sigma)$. Then there exist $(f_k)$ in the linear span of $\{P_n(\Sigma)\}$ such that $f_k \to f$ in $L^2(\Sigma)$, hence in $L^p(\Sigma)$. This shows the density for $p<2$. If $p>2$, let $h \in L^p(\Sigma)$ such that
$\int_\Sigma h P_n(\Sigma)\, d\sigma=0$ for every $n$. By the density of $span\, \{P_n(\Sigma)\}$ in $L^{p'}(\Sigma)$ then $\int_\Sigma h f(\Sigma)\, d\sigma=0$ for every $f \in L^{p'}(\Sigma)$, hence $h=0$.

\qed

We denote by $\sigma_p(\Sigma)$, $\{P_n(\Sigma)\}$ and $\{\lambda (P_n)\}$ the spectrum, the ($L^2$ normalized) eigenfunctions,  and the eigenvalues, listed according to their multiplicities, of $(-\Delta_0,D_p(\Delta_0, \Sigma))$, respectively.

\begin{lem} \label{eigen-S}
\begin{itemize}
\item[(i)]
$$\sigma_p(S_{N-1})=\{\lambda_n=n(n+N-2):\ n\in\N_0\}.$$
\item[(ii)]$$\sigma_p(S^+_{N-1})=\{\lambda_n=n(n+N-2):\ n\in\N\}$$
\end{itemize}
\end{lem}
{\sc Proof.} The first assertion is classical and the proof can be found in \cite[Chapter IX, Section 5.1]{vilenkin}. For the second assertion we refer to \cite[Proposition 4.5]{caldiroli}.
\qed

\begin{defi} \label{Fjp}
Let us fix $\Sigma \subset S_{N-1}$. For a given $J\subseteq \N_0$ we define 
\begin{equation*}  
F_{J,p}=F_{J,p}(\Sigma) =\overline{span}\{P_n(\Sigma):\ n\in J\}
\end{equation*}
where the closure is taken in $L^p(\Sigma)$ when $1\leq p<\infty$ and in $C_0(\Sigma)$, repectively.
\end{defi}
It is clear that 
$F_{J,p}$ is $\Delta_0$-invariant and the domain of 
${\Delta_0}_{|F_{J,p}}$ is given by $D_{p}(\Delta_0, \Sigma)\cap F_{J,p}$. We omit the label $\Sigma$ to shorten the notation, when no confusion may arise.
The following lemma is elementary.

\begin{lem}  \label{delta0fjp}
Let $1\leq p\leq \infty$. Then ${\Delta_0}_{|F_{J,p}}$ generates the analytic semigroup $(T_\Sigma(t)_{|F_{J,p}})_{t \ge 0}$  in $F_{J,p}$. Moreover
\begin{equation}  \label{restric}
\sigma_p({-\Delta_0}_{|F_{J,p}})=\{\lambda (P_n):\ n\in J\},
\end{equation}
where $\lambda (P_n)$ is the eigenvalue whose eigenfunction is $P_n(\Sigma)$.
\end{lem}
{\sc Proof. } Only the statement concerning the spectrum requires a proof. Since $P_n(\Sigma) \in F_{J,p}$, then $\lambda (P_n) \in \sigma_p({-\Delta_0}_{|F_{J,p}})$. For the converse note that the spectrum of ${-\Delta_0}_{|F_{J,p}}$ consists of eigenvalues, since the resolvent is compact, hence if $\lambda \in \sigma_p({-\Delta_0}_{|F_{J,p}})$, then $\lambda=\lambda(P_n)$ for some $n_0 \in \N_0$. If $n_0 \not \in J$, then $\int_{\Sigma}P_{n_0}(\Sigma)P_n(\Sigma)\, d\sigma=0$ for every $n \in J$ and then $P_{n_0}(\Sigma) \not \in F_{J,p}$ (the inner product is continuous with respect to the $L^p$-topology since  $P_{n_0}(\Sigma)$ is bounded). Then $n_0 \in J$ and the proof is complete.
\qed

Note that, since each eigenvalue can have more than one eigenfunction, different set of indeces leads to different spaces but not necessarily to different spectra.

The asymptotic behavior of $(T_\Sigma(t)_{|F_{J,p}})_{t \ge 0}$  in $F_{J,p}$ is determined by the first eigenvalue.
In the next lemma we assume that the numbers $\lambda(P_n)$ are listed in the increasing order. 
\begin{lem} \label {asymptotic}
Let $n$ be the smallest integer in $J$. There exists $M$ (depending on $n$ but not on $p$) such that for every $1 \le p \le \infty$
\begin{equation} \label{expdecay}
\|T_\Sigma(t)_{|F_{J,p}}\|_p \le M^{\big|1-\frac{2}{p}\big|}e^{-\lambda(P_n)\, t}.
\end{equation}
\end{lem}
{\sc Proof. } We may assume that $n>0$ and that $J=\{n, n+1,\dots \}$. Let $P$ be the $L^2$ orthogonal projection onto the linear span of $P_0(\Sigma), \dots , P_{n-1}(\Sigma)$ and observe that $P$ is bounded in $L^p(\Sigma)$ since the eigenfunctions are continuous. Then $Q=I-P$ is a bounded projection from $L^p(\Sigma)$ onto $F_{J,p}$ for every $1 \le p \le \infty$ (see also the proof of Lemma \ref{projection}). Then (\ref{expdecay}) follows if we prove that
\begin{equation} \label{expdecay1}
\|QT_\Sigma(t)\|_p \le M^{\big|1-\frac{2}{p}\big|}e^{-\lambda(P_n)\, t}.
\end{equation}
The above estimate holds (with equality) if $p=2$ since $Q$ has norm 1 and $\Delta_0$ is self-adjoint with eigenfunctions and eigenvalues $P_n(\Sigma)$, $\lambda (P_n)$, respectively.
By the Riesz-Thorin theorem it is sufficient therefore to prove (\ref{expdecay1}) for $p=1, \infty$.
Let $p=1$, consider $(T_\Sigma(t))_{t \ge 0}$ restricted on $F_{J,1}$ and let $S$ be the $L^2$ orthogonal projection on $Ker(\lambda(P_n)+\Delta_0)$, which is bounded in $F_{J,1}$ by the argument above.
Then 
$$F_{J,1}=Ker(\lambda(P_n)+\Delta_0)\oplus F_{K,1}$$
where $K=\{m+1,m+2, \dots \}$ for some $m \ge n$. Then $T_\Sigma (t)u=e^{-\lambda(P_n)\, t}u$ if $u \in Ker(\lambda(P_n)+\Delta_0)$ and $\|T_\Sigma (t)\| \le e^{-(\lambda(P_n)+\delta)\,t}$ for some $\delta >0$ on $F_{K,1}$, since $(T_\Sigma(t))_{t \ge 0}$ is analytic and its growth bound on $F_{K,1}$ coincides with the spectral bound which is strictly greater than $\lambda(P_n)$, by the preceeding lemma. Therefore $e^{\lambda(P_n)\, t} T_\Sigma (t) \to S$ in norm, as $t \to \infty$ and this shows (\ref{expdecay1}) for $p=1$. The proof for $p=\infty$ is the same.
\qed
Note that $M=1$ when $n=0$.

\subsection{The spaces $L^p_J({\cal C}_\Sigma)$}
If $X,Y$ are function spaces over $G_1, G_2$ we denote by $X\otimes Y$ the algebraic tensor product of $X,Y$, that is the set of all functions $u(x,y)=\sum_{i=1}^n f_i(x)g_i(y)$ where $f_i \in X, g_i \in Y$ and $x \in G_1, y\in G_2$.
If $T,S$ are linear operators  on $X,Y$ we denote by $T\otimes S$ the operator on $X\otimes Y$ defined by 
$$
T\otimes S \left (\sum_{i=1}^n f_i(x)g_i(y)\right)=\sum_{i=1}^n T f_i(x)Sg_i(y).
$$

\begin{defi} 
$L^p_{J}({\cal C}_\Sigma)$ $(1\leq p<\infty)$ and $C^0_J({\cal C}_\Sigma)$ $(p=\infty)$ are the closure of 
 $L^p(]0,+\infty[, \rho^{N-1}d\rho)\otimes F_{J,p}$ and $C_0^0(]0,+\infty[)\otimes F_{J,\infty}$, in $L^p(\Sigma)$ and $C^0(\Sigma)$, respectively. 
\end{defi}
We recall that the spaces $F_{J,p}$ have bee introduced in Definition \ref{Fjp}. 
In particular, if $J=\N_0$, then $L^p_{J}({\cal C}_\Sigma)=L^p({\cal C}_\Sigma)$ and $C^0_J({\cal C}_\Sigma)=C^0({\cal C}_\Sigma)$. To unify the notation we use $L^\infty_{J}({\cal C}_\Sigma)$ for $C^0_J({\cal C}_\Sigma)$.

\smallskip

The next lemma clarifies the structure of the spaces  $L^p_{J}({\cal C}_\Sigma)$ in some cases of interest.

\begin{lem} \label{projection}
Assume that the $L^2$ orthogonal projection $P: L^2(\Sigma) \to F_{J,2}$ extends to a bounded projection $P$ in $L^p(\Sigma)$. Then 
\begin{equation} \label{complement}
L^p({\cal C}_\Sigma)= L^p_{J}({\cal C}_\Sigma) \oplus L^p_{\N_0\setminus J}({\cal C}_\Sigma).
\end{equation}
In particular 
\begin{equation} \label{caratterizzazione}
 L^p_{J}({\cal C}_\Sigma)=\left \{u \in L^p({\cal C}_\Sigma): \int_\Sigma u(\rho\, \omega)P_j(\omega) \, d\sigma (\omega)=0\ {\rm for} \ \rho>0\ {\rm  and}\  j \not \in J\right \}.
\end{equation}
When $J$ is finite  
\begin{equation} \label{Jfinito1}
L^p_J({\cal C}_\Sigma)=\Bigl \{ u=\sum_{j \in J}f_j(\rho)P_j(\omega): f_j \in L^p(]0,+\infty[, \rho^{N-1}d\rho)\Bigr \}
\end{equation} and 
the projection $I\otimes P :L^p({\cal C}_\Sigma) \to L^p_{J}({\cal C}_\Sigma)$ is given by
\begin{equation} \label{Jfinito}
(I\otimes P) u=\sum_{j \in J}T_j u (\rho)\, P_j(\omega)=\sum_{j \in J}\left (\int_\Sigma u(\rho\, \omega)P_j(\omega) \, d\sigma (\omega) \right ) P_j(\omega).
\end{equation}
\end{lem}
{\sc Proof. } For $j \in \N_0$, $u \in L^p({\cal C}_\Sigma)$ let 
$$T_j u (\rho)=\int_\Sigma u(\rho\, \omega)P_j(\omega) \, d\sigma (\omega).$$
By H\"older inequality $\|T_ju\|_p \le c_j\|u\|_p$ (recall that $P_i \in L^\infty$). Then the right hand side of (\ref{caratterizzazione}) is $\cap_{j \not \in J}Ker\, (T_j)$. Since $L^p_{J}({\cal C}_\Sigma)$ is the closure of functions of the form $\sum_{j \in J}f_j(\rho)P_j(\omega)$, it follows that $L^p_{J}({\cal C}_\Sigma) \subset \cap_{j \not \in J}Ker\, (T_j)$. Similarly $L^p_{\N_0 \setminus J}({\cal C}_\Sigma) \subset \cap_{j \in J}Ker\, (T_j)$ and hence 
$$
L^p_{J}({\cal C}_\Sigma) \cap L^p_{\N_0 \setminus J}({\cal C}_\Sigma)  \subset \cap_{j \in N_0}Ker\, (T_j)=\{0\},
$$
since the linear span of the functions $\{P_j,\,  j \in\N_0\}$ is dense in $L^{p'}(\Sigma)$ by Lemma \ref{eigen-S} (and $w^*$-dense if $p'=\infty$).

 Let $Q: L^2(\Sigma) \to F_{\N_0\setminus J,2}$ be the orthogonal projection. Since $P+Q$ is the identity on $L^2( \Sigma)$ and $P$ is bounded in  $L^p( \Sigma)$, then $Q$ is bounded in  $L^p( \Sigma)$, too.
Next note that $I\otimes P, I\otimes Q$ are bounded projections in $L^p({\cal C}_\Sigma)$ and that
$$I\otimes P \left (\sum_j f_j(\rho)P_j(\omega)\right )= \sum_{j \in J}f_i(\rho)P_i(\omega),
$$
where the sums are finite, and similarly for $I\otimes Q$. It follows that the ranges of $I\otimes P, I\otimes Q$ are 
$ L^p_{J}({\cal C}_\Sigma)$, $ L^p_{\N_0\setminus J}({\cal C}_\Sigma)$ respectively and that $I\otimes P+I\otimes Q$ is the identity of $L^p({\cal C}_\Sigma)$. Then (\ref{complement}) follows.
Equality (\ref{caratterizzazione}) follows from (\ref{complement}) since $L^p_{J}({\cal C}_\Sigma) \subset \cap_{j \not \in J}Ker\, (T_j)$, $L^p_{\N_0 \setminus J}({\cal C}_\Sigma) \subset \cap_{j \in J}Ker\, (T_j)$ and the subspaces  $\cap_{j \not \in J}Ker\, (T_j)$,  $\cap_{j \in J}Ker\, (T_j)$ intersect only at 0.

\smallskip

Finally, we assume that $J$ is finite and observe that  (\ref{Jfinito1}) holds if and only if its right hand side is closed. Let  $u^k=\sum_{j \in J}f^k_j(\rho)P_j(\omega)$ converge to $u$ in $L^p(\R^N)$. Then $f_j^k=T_j(u^k)$ converges to some $f_j$ in  $ L^p(]0,+\infty[, \rho^{N-1}d\rho)$ and therefore $u=\sum_{j \in J}f_j(\rho)P_j(\omega)$. This proves (\ref{Jfinito1}). Identity (\ref{Jfinito}) holds since both sides are continuous in $L^p({\cal C}_\Sigma)$, $J$ being finite, and coincide on finite sums  $\sum_i f_i(\rho)P_i(\omega)$.
\qed

\begin{os} {\rm Note that the inclusion  
$$ L^p_{J}({\cal C}_\Sigma) \subset \left \{u \in L^p({\cal C}_\Sigma): \int_\Sigma u(\rho\, \omega)P_j(\omega) \, d\sigma (\omega)=0\ {\rm for} \ \rho>0\ {\rm  and}\  j \not \in J\right \}
$$
holds without assuming the boundedness of the projection $P$.}
\end{os}

Special situations of interest where Lemma \ref{projection} applies are the following.

\begin{defi} \label{Lp<n} {\rm We fix $n \in \N_0$ and define $L^p_{\le n}, L^p_n, L^p_{>n}$ as the closure in $L^p(\R^N)$ of the functions
$
u=\sum_j f_j(\rho) P_j(\omega)
$
where the sums are finite, $f_j\in L^p(]0,\infty[, \rho_{N-1}\, d\rho)$ and $(P_j)$ constitute a basis for spherical harmonics of order $\le n$, $n$ and $> n$ respectively. The spaces $L^p_{<n}, L^p_{\ge n}, L^p_{\neq n}$ are defined similarly. } \\
\end{defi}

Note that $L^p_0$ consists of radial functions and that $L^p(\R^N)=L^p_{\le n}\oplus L^p_{>n}$.

\subsection{The spectrum of $A$}

The following result follows from well- known and elementary facts, see \cite[AI, Section 3.7]{nagel}.

\begin{prop} \label{analyt}
For $1\le p \le \infty$, let  $D_p(\Gamma)$ and $D_p({\Delta_0}_{|F_{J,p}})$ be the domains of $\Gamma$ and ${\Delta_0}_{|F_{J,p}}$ introduced in the previous subsection. Then the closure of the operator $(A,D_p(\Gamma)\otimes D_p({\Delta_0}_{|F_{J,p}}))$  generates a strongly continuous  analytic semigroup $(T_{p,J,\Sigma}(t))_{t \ge 0}$
in $L^p_J({\cal C}_\Sigma)$.  Let $n$ be the smallest integer in $J$. There exists $M$ (depending on $n$ but not on $p$) such that for every $1 \le p \le \infty$
\begin{equation} \label{expdecay2}
\|T_{p,J,\Sigma}(t)\|_p \le M^{\big|1-\frac{2}{p}\big|}e^{-(\omega_p+\lambda(P_n))\, t},
\end{equation}
where $\omega_p$ is defined in (\ref{defomegap}).
\end{prop}
{\sc Proof.} Observe first that 

$$A=\Gamma\otimes Id +Id\otimes{\Delta_0}_{|F_{J,p}}$$ on $D_p(\Gamma)\otimes D_p({\Delta_0}_{|F_{J,p}})$.
Let $(S(z))_{z \in \Omega}$ and $(T_J(z))_{z\in \Omega}$, where $\Omega$ is a suitable sector in the complex plane,  be  the analytic semigroups generated respectively by $\Gamma$ in $L^p(]0,+\infty[, \rho^{N-1}d\rho)$ and ${\Delta_0}_{|F_{J,p}}$ in $F_{J,p}$. The family $(S(z)\otimes T_J(z))_{z\in \Omega}$ extends to a strongly continuous analytic semigroup $(T_{p,J,\Sigma}(t))_{t \ge 0}$ on $L^p_J({\cal C}_\Sigma)$.
Moreover the generator of $(S(t)\otimes T_J(t))_{t\geq 0}$ is given by  the closure of the operator $$\Gamma\otimes Id +Id\otimes{\Delta_0}_{|F_{J,p}}$$ defined on  the core $D_p(\Gamma)\otimes D_p({\Delta_0}_{|F_{J,p}})$.
Finally, since by Proposition \ref{tensornorm} 
$$
\|T_{p,J,\Sigma}(t)\|_p=\|S(t)\|_p\|T_J(t)\|_p,
$$
 (\ref{expdecay2}) follows from Proposition \ref{spec-rad} and Lemma \ref{asymptotic}.
\qed

Note that $M=1$ when $n=0$.
We denote by $A_{p,J,\Sigma}$ the closure of  $(A,D_p(\Gamma) \otimes D_p({\Delta_0}_{|F_{J,p}}))$ 
in $L^p_{J}({\cal C}_\Sigma)$. When $I=\N_0$ we write $A_{p,\Sigma}$ for $A_{p,J,\Sigma}$ and $T_{p,\Sigma}(t)$ for $T_{p,J,\Sigma}(t)$.

\begin{cor} \label{restriction}
$T_{p,J,\Sigma} (t)$ is the restriction of $T_{p,\Sigma}(t)$ to $L^p_J({\cal C}_\Sigma)$ and its generator $A_{p,J,\Sigma}$ is the part of $A_{p,\Sigma}$ in $L^p_J({\cal C}_\Sigma)$.
\end{cor}
{\sc Proof. } Keeping the notation of the proof of Proposition \ref{analyt}, if $(T_\Sigma (t))_{t \ge 0}$ is the semigroup generated by $\Delta_0$ in $L^p(\Sigma)$, then $T_J(t)=(T_\Sigma (t))_{|F_{J,p}}$, by Lemma \ref{delta0fjp}. Then 
the restriction of $S(t) \otimes T_\Sigma (t)$ on $D_p(\Gamma)\otimes D_p({\Delta_0}_{|F_{J,p}})$ coincides with $S(t) \otimes T_J(t)$ and hence $T_{p,J,\Sigma} (t)$ is the restriction of $T_{p,\Sigma}(t)$ to $L^p_J({\cal C}_\Sigma)$. The second statement follows from basic semigroup theory.
\qed

In the next proposition we show that smooth functions are a core for $A_{p,J,\Sigma}$.

\begin{prop}  \label{core}
The set $$\{u\in C_c^\infty({\ov{\cal C}_\Sigma}\setminus\{0\}):\quad u\equiv 0\ \textrm{on}\ \partial {\cal C}_\Sigma\}$$ is a core for $A_{p,J,\Sigma}$.
\end{prop}
{\sc Proof.} Observe that,  since $C_c^\infty(]0,+\infty[)$  is dense  in $D_p(S^{-1}\Gamma S)$ (see Proposition \ref{spec-rad}), then $C_c^\infty(]0,+\infty[)$ is  also dense in $D_p(\Gamma)$. Moreover $\{u\in C^\infty(\Sigma)|\ u\equiv 0 \ \textrm{on}\ \partial \Sigma\}$ is dense in $ D_p({\Delta_0}_{|F_{J,p}})$. By Proposition \ref{analyt}, $D_p(\Gamma)\otimes  D_p({\Delta_0}_{|F_{J,p}})$ is a core for $A_{p,J,\Sigma}$. It follows that 
$$C_c^\infty(]0,+\infty[)\otimes \{u\in C^\infty(\Sigma)|\ u\equiv 0 \ \textrm{on}\ \partial \Sigma\}$$
is dense in $D_p(A)$. Observing that $$C_c^\infty(]0,+\infty[)\otimes \{u\in C^\infty(\Sigma)|\ u\equiv 0 \ \textrm{on}\ \partial \Sigma\}\subseteq \{u\in C_c^\infty({\ov{\cal C}_\Sigma}\setminus\{0\}):\quad u\equiv 0\ \textrm{on}\ \partial {\cal C}_\Sigma\}$$ we get the claim. \qed
We can now prove the main result of this section.

\begin{teo} \label{main} 
Let $1\leq p\leq \infty$. Then
$$\sigma(A_{p,J,\Sigma})=\cup_{n\in J}({\cal P}_p-\lambda (P_n)),$$
where ${\cal P}_p$ is the parabola defined in (\ref{spettrogamma}).
\end{teo}
{\sc Proof.} Let $\lambda \not \in \cup_{n\in J}({\cal P}_p-\lambda (P_n))$ and fix $n \in N_0$ such that
\begin{equation} \label{boundgamma}
-\omega_p-\lambda(P_k) < Re\, \lambda \quad {\rm for\  every\ } k > n.
\end{equation} 
According to Lemma \ref{projection} we write
$L^p_J({\cal C}_\Sigma)= L^p_{J_n}({\cal C}_\Sigma) \oplus L^p_{J\setminus J_n}({\cal C}_\Sigma)$,
with $J_n=J\cap \{0,1,\dots ,n\}$. Since both  $L^p_{J_n}({\cal C}_\Sigma)$  and  $L^p_{J\setminus J_n}({\cal C}_\Sigma)$ are $A_{p,J, \Sigma}$ invariant, then $\lambda \in \rho(A_{p,J,\Sigma})$ if and only if  $\lambda \in \rho(A_{p,J_n,\Sigma})$ and  $\lambda \in \rho(A_{p,J \setminus J_n,\Sigma})$. The second inclusion follows immediately from (\ref{expdecay2}) with $J\setminus J_n$ instead of $J$, since $Re\, \lambda$ is greater than the growth bound of $(T_{p, J \setminus J_n, \Sigma})_{t \ge 0}$, by  (\ref{boundgamma}). Concerning the first inclusion we note that
$$
 L^p_{J_n}({\cal C}_\Sigma) =\oplus_{i=0}^n  L^p_{J_i}({\cal C}_\Sigma) 
$$
where $J_i=J\cap \{i\}$ and that each $ L^p_{J_i}({\cal C}_\Sigma) $ is $A_{p,J,\Sigma}$ invariant. Moreover, $\lambda-A_{p,J, \Sigma}$ coincides with $\left (\lambda+\lambda (P_i)-\Gamma\right ) \oplus I$ on $ L^p_{J_i}({\cal C}_\Sigma) $, hence it is invertible on it, since $\lambda+\lambda(P_i) \not \in {\cal P}_p$ by assumption.
This shows that $\lambda \in \rho (A_{p,J, \Sigma})$, hence
$$\sigma(A_{p,j,\Sigma})\subseteq \sigma_p(\Gamma)+\sigma_p({\Delta_0}_{|F_{J,p}})=\cup_{n\in J}({\cal P}_p-\lambda (P_n)).$$
Let us now prove the reverse inclusion, assuming first that $1 \le p <\infty$. Fix $n\in J$,$\xi \in \R$  let $\lambda_n=\lambda (P_n)$ and set $\mu= -\xi^2+i\xi(N-2+c-\frac{2N}{p})-\omega_p-\lambda_n:=\lambda-\lambda_n$. Set also $u(x)=\rho^{-\frac{N}{p}+i\xi}P(\omega)$ where $P$ is an eigenfunction of ${-\Delta_0}_{|F_{J,p}}$ corresponding to the eigenvalue $\lambda_n$. 
Then
$$A_{p,J,\Sigma}\,  u=(\Gamma\otimes Id +Id\otimes{\Delta_0}_{|F_{J,p}} )u=\lambda u$$
but $u$ does not belong to $L^p(\R^N)$. We  approximate $u$ by a sequence $(u_k)_{k\in\N} \subset D_p(\Gamma)\otimes D_p({\Delta_0}_{|F_{J,p}})$ such that $\lambda u_k-Lu_k \to 0$, $u_k \not \to 0$ as $k$ goes to infinity and we deduce that $\lambda\in \sigma(A_{p,J,\Sigma})$. 
Let $\phi_k(\rho)\in C_c^\infty(\R)$, $0\leq \phi_k\leq $ such that $\phi_k(\rho)=0$ if $r\le\frac{1}{2k}$, $\phi_k(\rho)=1$ if $\frac{1}{k}\leq \rho\leq k$, $\phi_k(\rho)=0$ if $\rho\geq 2k$ and
$$|\phi_k'|\leq C\left(k\chi_{\left[\frac{1}{2k},k\right]}+\frac{1}{k}\chi_{[k,2k]}\right),$$
$$|\phi_k''|\leq C\left(k^2\chi_{\left[\frac{1}{2k},k\right]}+\frac{1}{k^2}\chi_{[k,2k]}\right)$$
and set $u_k(x)=\phi_k(\rho)\rho^{-\frac{N}{p}+i\xi}P(\omega)$. We have that
\begin{equation} \label{Norma}
\int_{{\cal C}_\Sigma}|u_k|^pdx\geq\int_{\Sigma}|P(\omega)|^pd\sigma\int_{\frac{1}{k}}^k\frac{1}{\rho}d\rho\geq C \log k.
\end{equation}
Moreover 
\begin{align*}
Au_k=A(\phi_k u)=\lambda u_k+ |x|^2(\Delta\phi_k) u+2|x|^2\nabla\phi_k\cdot\nabla u+cx\cdot \nabla \phi_k u.
\end{align*}
By the definition of $u$,
$$\nabla\phi_k\cdot\nabla u=\phi'_k(\rho)\frac{x}{\rho}\cdot\nabla u=\phi'_k(\rho)\frac{\partial u}{\partial \rho}=\phi'_k(\rho)\left(-\frac{N}{p}+i\xi\right)\rho^{-\frac{N}{p}+i\xi-1}P(\omega).$$
Therefore 
\begin{align*}
 \|\rho^2\nabla\phi_k\cdot\nabla u\|_p^p&\leq C \left(k^p\int_\frac{1}{2k}^\frac{1}{k}\rho^{2p-N-p+N-1}d\rho+ \frac{1}{k^p}\int_k^{2k} \rho^{2p-N-p+N-1}d\rho\right)\\
&= C \left(k^p\int_\frac{1}{2k}^\frac{1}{k}\rho^{p-1}dr+ \frac{1}{k^p}\int_k^{2k} \rho^{p-1}d\rho\right)\leq C.
\end{align*}
Similarly
\begin{align*}
 \|r^2\Delta\phi_k u\|_p^p&\leq C \left(k^{2p}\int_\frac{1}{2k}^\frac{1}{k}\rho^{2p-N+N-1}d\rho+ \frac{1}{k^{2p}}\int_k^{2k} \rho^{2p-N+N-1}d\rho\right)\\
&= C \left(k^{2p}\int_\frac{1}{2k}^\frac{1}{k}\rho^{2p-1}dr+ \frac{1}{k^{2p}}\int_k^{2k} \rho^{2p-1}d\rho\right)\leq C
\end{align*}
and, since $x\cdot \nabla \phi_k u=\rho\phi'_ku$, we have also $\|x\cdot \nabla \phi_k u\|_p^p\leq C$.
It follows that $\|Au_k-\lambda u_k\|_p\leq C$ for some positive constant $C$. Then, setting $v_k=\frac{u_k}{\|u_k\|_p}$, by (\ref{Norma}) and the computation above  $\|Av_k-\lambda v_k\|_p$ goes to $0$ as $k$ goes to infinity and this implies that $\lambda\in\sigma(A_{p,J,\Sigma})$.
The proof for $p=\infty$ is similar.
\qed
 Observe that $\sigma(A_{p,J,\Sigma})$ is real if and only if $N-2+c-2N/p=0$ hence in the self-adjoint case (that is when $c=2$) if and only if $p=2$.

\begin{os} {\rm The inclusion 
$$\sigma(A_{p,j,\Sigma})\subseteq \sigma_p(\Gamma)+\sigma_p({\Delta_0}_{|F_{J,p}})=\cup_{n\in J}({\cal P}_p-\lambda (P_n))$$
follows also from the more general result  \cite[Theorem 7.3]{arendt1} since the semigroups generated by $\Gamma$ and ${\Delta_0}_{|F_{J,p}}$ are analytic and commute.
 }
\end{os}

\section{The operator $A$ in the whole space and in the half-space}
In this section we complete the analysis of $A$ by computing the growth bound in $\R^N$ and in $\R_+^N$, that is when $\Sigma=S_{N-1}, S_{N-1}^+$, respectively, and describing the domain in the case of the whole space.

\subsection{Part I: $A$ in $\R^N$}
Here $\Sigma=S_{N-1}$, $J=\N_0$ and we wrire $A_p$ for $A_{p,J,\Sigma}$.
The following result is a particular case of Theorem \ref{main}

\begin{prop} \label{specRN}
Let $1\leq p\leq \infty$. Then
$$\sigma(A_p)=\cup_{n\in\N_0}({\cal P}_p-\lambda_n)$$
where $\lambda_n$ and ${\cal P}_p$ are as in Lemma \ref{eigen-S} (i) and Proposition \ref{spec-rad}, respectively.
\end{prop}

The operator $A_p$ generates a positive semigroup $\Tt$ in $L^p(\R^N)$ which is independent of $p$. Moreover $C_c^\infty (\R^N\setminus \{0\})$ is a core for $A_p$.  These results follow from Propositions \ref{analyt} and \ref{core}. 
As in \cite[Section 3, Section 4]{met-spi1}, see also \cite{for-lor}, one can define for $1<p<\infty$
 \begin{align*}
D_p(A)=& \{u \in L^p(\R^N)\cap W^{2,p}\left (\R^N \setminus B_\eps\right)\ {\rm for\ every\ } \eps>0: 
 |x|\nabla u, |x|^2 D^2u \in L^p(\R^N)\},
\end{align*}
and show that it coincides with the maximal one
\begin{align*}
D_{p,max}(A)=& \{u \in L^p(\R^N)\cap W^{2,p}\left (\R^N \setminus B_\eps\right)\ {\rm for\ every\ } \eps>0: Au \in L^p(\R^N)\}.
\end{align*}
Moreover $(A,D_p(A))$ generates an analytic semigroup which coincides with that of Proposition \ref{analyt}, as shown in the next Proposition.
The domain in $C_0^0(\R^N)$, the spaces of all continuous functions vanishing at $0,\infty$, coincides with the maximal one

\begin{equation*}
D^0_{max}(A)= \{u \in C_0^0(\R^N)\cap W_{loc}^{2,p}(\R^N \setminus \{0\} {\rm \ for \ every\ }  p<\infty: Au \in C_0^0(\R^N)\}.
\end{equation*}

\begin{prop} \label{domini}
If $1<p<\infty$, then  
the closure of $(A,D_p(\Gamma) \otimes D_p(\Delta_0))$ coincides with $(A,D_p(A))$. If $p=\infty$ the corresponding closure coincides with $(A,D_{max}^0(A))$.
\end{prop}
{\sc Proof.} Let $1<p<\infty$.
Since $D_p(A)=D_{p,max}(A)$, it follows that $D_p(\Gamma)\otimes D_p(\Delta_0)\subseteq D_p(A)$. Hence the closure of 
$(A,{D_p(\Gamma)\otimes D_p(\Delta_0)})$ is contained in $(A,D_p(A))$.
Since both operators generate a semigroup in $L^p(\R^N)$, the equality follows.
The proof in $C_0^0(\R^N)$ is identical.
\qed

\bigskip
\noindent Next, we   estimate the growth bound of $\Tt$.
Even though the result below can be deduced from Propostion \ref{analyt}, see equation (\ref{expdecay2}), we prefer to give a direct proof which shows the equivalence with Hardy inequalities and which applies also in the case of the half-space, see Proposition \ref{dissipativity2}.

\begin{prop} \label{dissipativity1}
Set $\omega_p=\frac{N}{p^2}\left[p(N-2+c)-N\right]$, see (\ref{defomegap}).  Then $\|T(t)\|_p \le e^{-\omega_p t}$. The constant $\omega_p$ is sharp. 
\end{prop}
{\sc Proof.}  Consider first the case $1< p<\infty$.
Since $C_c^\infty$ is a core for $A_p$, by Proposition \ref{core},
it suffices to show the dissipativity estimate 
$$-\int_{\R^N}Au |u|^{p-2}u\;
dx\geq \omega_p \int_{\R^N}|u|^p\;
dx$$ for every $u\in C_c^\infty$.
Setting
$u^\star=u|u|^{p-2}$ we multiply  $Lu$  by
$u^\star$ and integrate over $\R^N$. The integration by parts is
straightforward when $p\geq 2$. For $1<p<2$, $|u|^{p-2}$
becomes singular near the zeros of $u$ but integrating by parts is still allowed, see
\cite{met-spi}. We get
\begin{align*}
-\int_{\R^N}&Au\, u^\star\;
dx=(p-1)\int_{\R^N}|x|^2|u|^{p-2}|\nabla
u|^2\;dx+ (2-c)\int_{\R^N}\nabla u |u|^{p-2}u\cdot x\; dx\\&=(p-1)\int_{\R^N}|x|^2|u|^{p-2}|\nabla
u|^2\;dx+ \left(\frac{2-c}{p}\right)\int_{\R^N}\nabla |u|^p\cdot x\; dx\\&=(p-1)\int_{\R^N}|x|^2|u|^{p-2}|\nabla
u|^2\;dx- N\left(\frac{2-c}{p}\right)\int_{\R^N}|u|^p\; dx.
\end{align*}
By Hardy inequality (\ref{hardy2}) with $\beta=2$,
\begin{align*}
-\int_{\R^N}&Au\, u^\star\;
dx\geq \left[(p-1)\frac{N^2}{p^2}- N\left(\frac{2-c}{p}\right)\right]\int_{\R^N}|u|^p\; dx=\omega_p \int_{\R^N}|u|^p\; dx
\end{align*}
and therefore
$$-\int_{\R^N}Au |u|^{p-2}u 
dx\geq \omega_p \int_{\R^N}|u|^p\;
dx.$$
Observe that all inequalities above are equalities, except for Hardy inequality. Hence $\omega_p$ is sharp since the constant in  (\ref{hardy2}) is sharp.

By standard semigroup theory the above estimate is equivalent to $\|T(t)f\|_p \le e^{-\omega_p\, t}\|f\|_p$ for every $f \in L^p$. Letting $p \to 1$ we get the same estimate in $L^1$ (with $\omega_1$ instead of $\omega_p$).

If $p=\infty$ (that is in $C_0^0(\R^N)$) the result follows from the maximum principle and the positivity of $\Tt$, see Appendix B.
\qed

\begin{os} \label{L1}
{\rm The estimate  $\|T(t)f\|_1 \le e^{-\omega_1\, t}\|f\|_1$ shows that $\Tt$ extends to a semigroup in $L^1$. The strong continuity follows from \cite[Proposition 4]{Voigt}.}

\end{os}

\begin{os} {\rm
\begin{itemize}
 \item[(i)] $\omega_\infty=0$ and $\omega_1=(c-2)N$;
\item[(ii)] $\omega_p\geq 0$ iff $p\geq \frac{N}{N-2+c}$. Moreover $\omega_p$ attaints its maximum value at $\ov{p}=\frac{2N}{N-2+c}$ and $\omega_{\ov{p}}=\left(\frac{N-2+c}{2}\right)^2$;
\item [(iii)] if $\frac{N}{N-2+c}<p<\infty$, then $\omega_p>0$ and $A_p$ is invertible in $L^p$. 
\end{itemize}}
 \end{os}

\subsection{Part II: $A$ in $\R^N_+$}
In this case $\Sigma=S_{N-1}^+$ and we write $A_p^+$ for $A_{p,S_{N-1}^+}$. The following result is a particular case of Theorem \ref{main} but differs from Proposition \ref{specRN}. In particular the spectral bound of $A_p$ is $-\omega_p$ whereas the spectral bound of $A_p^+$ is $-\omega_p-\lambda_1=-\omega_p-(N-1)$.

\begin{prop} \label{specRN+}
Let $1\leq p\leq \infty$. Then
$$\sigma(A^+_p)=\cup_{n\in\N}({\cal P}_p-\lambda_n)$$
where $\lambda_n$ and ${\cal P}_p$ are as in Lemma \ref{eigen-S} (ii) and Proposition \ref{spec-rad}, respectively.
\end{prop}

The operator $A_p^+$ generates a positive semigroup $(T^+(t))_{t \ge 0}$ in $L^p (\R^N_+)$  which is independent of $p$: this follows from Proposition \ref{analyt}. The semigroup $(T^+(t))_{t \ge 0}$ is pointwise dominated by (the restriction to $\R^N_+$ of) the semigroup $\Tt$ generated by $A_p$. However, its grouth bound is strictly smaller than $-\omega_p$, defined in Proposition \ref{dissipativity1}, since Hardy inequality in the half-space holds with a better constant than in the whole space, see Proposition \ref{hardy2} .

\begin{prop} \label{dissipativity2}
Set 
\begin{equation} \label{defomegap+}
\omega_p^+=\frac{N}{p^2}\left[p(N-2+c)-N\right]+\ds\frac{4(p-1)(N-1)}{p^2}=\omega_p+\ds\frac{4(p-1)(N-1)}{p^2}.
\end{equation}
Then $\|T^+(t)\|_p \le e^{-\omega_p^+ t}$.  The constant $\omega_p^+$ is sharp.
\end{prop}
{\sc Proof.} The statement is equivalent to the dissipativity estimate
$$-\int_{\R^N_+}A_p^+u|u|^{p-2}u\, dx\geq\omega_p^+\int_{\R^N_+}|u|^p\, dx$$ for every $u \in C_\infty (\overline{\R^N_+} \setminus \{0\})$ such that $u\equiv 0$ on $\R^{N-1}$, since the last space is a core for $A_p^+$, by Proposition \ref{core}. The proof is a repetition of that of Proposition \ref{dissipativity1}, using (\ref{WHardy2-half}) instead of (\ref{WHardy2-RN}).
As in the proof of Proposition \ref{dissipativity1}, the sharpness of $\omega_p^+$ follows from the sharpness of the constant in Hardy inequality (\ref{WHardy2-half}).
\qed

\begin{os} \label{noncontractivity}
{\rm Note  that $\omega_p^*=\omega_p$ if $p=1,\infty$. It is worth mentioning that the dissipativity constant $-\omega_p^+$ is greater that the spectral bound $s_p :=-\omega_p-(N-1)$ and coincides with it if and only if $p=2$. Since the semigroup  $(T^+(t))_{t \ge 0}$ is analytic, $s_p$ coincides with the growth bound. This means that for every $\eps>0$ there exists a constant $C_\eps>0$ such that
$\|T^+(t)\|_p \le C_\eps e^{(s_p+\eps)t}.
$
However the estimate $\|T^+(t)\|_p \le  e^{ \mu t}$ holds only if $\mu \ge -\omega_p^+$.
}

\end{os}

\section{Special cases and generalizations}

\subsection{Rellich inequalities for Schr\"odinger operators}
Rellich inequalities for Schr\"odinger operators can be deduced from Theorem \ref{Rellich-p} by settting $c=0$, so that
$$L=\Delta-\ds\frac{b}{|x|^2}.$$ For simplicity we assume that $b \in \R$ and $N \ge 2$ and we observe that if $b=0$ the operator $L$ reduces to the Laplacian.
Note that
\begin{equation*} 
\gamma_p(\alpha,0)=\Bigl(\frac{N}{p}-2+\alpha\Bigr)\Bigl(\frac{N}{p'}-\alpha\Bigr)
\end{equation*}
and 
\begin{equation*} 
{\cal P}_{p,\alpha,0}:=
\left\{\lambda=-\xi^2+i\xi\Bigl (N\Bigl (1-\frac2p \Bigr)+2-2\alpha\Bigr)-\gamma_p(\alpha)\;;\;\xi\in \R\right\}.
\end{equation*}

\begin{prop} \label{Rellich-Scrod}
Let $N\geq 2$, $\alpha, b\in\R$, $1 \le p \le \infty$.
Then there exists a positive constant $C=C(N,\alpha,  p, b)$ such that
\begin{equation*} 
\left\||x|^\alpha \left(\Delta u-\ds\frac{b}{|x|^2}u\right)\right\|_p \geq C\||x|^{\alpha-2} |u|\|_p
\end{equation*}
holds for every $u\in C_c^\infty(\R^N \setminus \{0\})$ if and only if $b+\lambda_n \not \in {\cal P}_{p,\alpha,0}$ for every $n \in \N_0$.
If $(N(1-2/p) +2-2\alpha \not =0$, then ${\cal P}_{p,\alpha,0}$ is a non degenerate parabola with vertex at $(-\gamma_p(\alpha,0),0)$ and the above condition reads $b+\gamma_p(\alpha,0)+\lambda_n \not = 0$ for every $n\in \N_0$ or, equivalently, since $\lambda_n=n^2+(N-2)n$
\begin{equation} \label{condschr}
b+\Bigl(n+\frac{N}{2}-1\Bigr)^2 \neq\Bigl(N\Bigl (\frac12-\frac1p\Bigr)+1-\alpha\Bigr)^2.
\end{equation}
 However, if $(N(1-2/p) +2-2\alpha  =  0$, then ${\cal P}_{p,\alpha,0}$  coincides with the semiaxis $]-\infty,-\gamma_p(\alpha,0)]$ and the condition becomes $b+\gamma_p(\alpha,0)>0$.\\
When  $b+\gamma_p(\alpha,0)>0$
the optimal constant is given by
$$C=b+\gamma_p(\alpha,0)=b+\Bigl(\frac{N}{p}-2+\alpha\Bigr)\Bigl(\frac{N}{p'}-\alpha\Bigr).$$
\end{prop}

It is worth-mentioning that Rellich inequalities hold for large negative $b$ satisfiyng (\ref{condschr})
 for every $n \in \N_0$ if the parabola ${\cal P}_{p,\alpha,0}$ is non degenerate. If $\alpha=0$ then ${\cal P}_{p,\alpha,0}$ degenerates if and only if $p_0=2N/(N+2)$ and Rellich inequalities hold if and only if $b+\gamma_{p_0}(0,0)=b+\left (\frac{N}{2} -1\right )^2>0$.

If $p=2$, the parabola ${\cal P}_{2,\alpha,0}$ degenerates into a half-lineif and only if $\alpha=1$, hence the inequality
\begin{equation*} 
\left\||x| \Delta u-\ds\frac{b}{|x|^2}u\right\|_2 \geq C\||x|^{-1} u\|_2
\end{equation*}
holds if and only if $b+(N/2-1)^2>0$.

When $p=2$ better results can be deduced from Section 2, concerning the computation of the best constants, as in \cite{caldiroli}. In the classical case $\alpha=0$, then $\gamma_2(0,0)=(N/4)(N-4)$ and Rellich inequalities 
\begin{equation*} 
\left\| \Delta u-\ds\frac{b}{|x|^2}u\right\|_2 \geq C\||x|^{-2} |u|\|_2
\end{equation*}
hold if and only if $b+\Bigl(n+\frac{N}{2}-1\Bigr)^2 \neq 1$ for every $n \in \N_0$.

\subsection{Rellich and Calder\'on-Zygmund inequalities in $L^p_J(\R^N)$}
The failure of Rellich and  Calder\'on-Zygmund inequalities  for some values of $\alpha$ (depending on $N,p,b,c$) is determined by certain subspaces defined by spherical harmonics of low order. Discarding these subspaces, Rellich and  Calder\'on-Zygmund inequalities continue to hold even though they fail in the whole $L^p$. This phenomenon holds also in the extreme cases $p=1,\infty$. 

Let $J \subset \N_0$ and $F_{J,p}$ and $L^p_J(\R^N)$ be defined as in Section \ref{spectrumA} (with $\Sigma=S_{N-1}$) by selecting for every $j\in J$ a (different) spherical harmonic $P_j$. We denote by $\lambda (P_j)$ the eigenvalue corresponding to $P_j$ and we suppose that they are listed in the increasing order.

By using the results on the spectrum of $A$  in $L^p_J(\R^N)$
(see Theorem \ref{main}, Section \ref{spectrumA}) and by arguing as in Theorems \ref{Rellich-p} and \ref{Cald-Zyg}, we  can improve Rellich and Calder\'on-Zygmund inequalities.
\begin{teo} \label{RellichFJ}
Let $N\geq 2$, $\alpha, b,\ c\in\R$, $1 \le p \le \infty$. Let $\gamma_p(\alpha,c)$ and the parabola ${\cal P}_{p,\alpha,c}$ be defined as in (\ref{gammap}) and (\ref{Pp}).
 Then there exists a positive constant $C=C(N,\alpha,  p, c, b)$ such that
\begin{equation*} 
\||x|^\alpha Lu\|_p \geq C\||x|^{\alpha-2} |u|\|_p
\end{equation*}
holds for every $u\in C_c^\infty(\R^N \setminus\{0\}) \cap L^p_J(\R^N)$, if and only if $b+\lambda (P_j) \not \in {\cal P}_{p,\alpha,c}$ for every $j\in J$. Moreover, if $n=\min J$ and $b+\lambda (P_n)+\gamma_p(\alpha,c)>0$ for every $1 \le p \le \infty$, the best constant $C$ above satisfies
\begin{equation} \label{bestLpJ}
c^{\big |1-\frac2p |}\left (b+\lambda (P_n)+\gamma_p(\alpha,c)\right ) \le C\le b+\lambda (P_n)+\gamma_p(\alpha,c)
\end{equation}
where $c>0$ depends on $n, N$ but not on $p$.
\end{teo}
{\sc Proof.} The proof is identical to that of Theorem \ref{Rellich-p} but this time the auxiliary operator  $\tilde L$ and its spectrum is considered in  $L^p_J(\R^N)$ instead of $L^p(\R^N)$. An application of Theorem \ref{main} instead of Proposition \ref{specRN} concludes the first part of the proof. Concerning the estimate of the best constant $C$, we recall that, as in Theorem \ref{Rellich-p}, $C^{-1}=\|(b-\tilde L)^{-1}\|_p$, where the operator $\tilde L$ is defined in (\ref{deftildeL}) and it is considered in $L^p_J(\R^N)$. Since, by Theorem \ref{main},
$$
\sigma (\tilde L)=\cup_{j \in J}({\cal P}_{p,\alpha,c}-\lambda(P_j))
$$
and $-\gamma_p(\alpha,c)$ is the vertex of the parabola ${\cal P}_{p,\alpha,c}$, the assumption $b+\lambda (P_n)+\gamma_p(\alpha,c)>0$ for every $1 \le p \le \infty$ implies the spectrum of $\tilde L$ lies in the half-plane $\{Re\, \lambda <b\}$, hence $b$ belongs to the resolvent set for every $1 \le p\le \infty$ and 
$C=(\|(b-\tilde L)^{-1}\|_p)^{-1} \le dist (b,\sigma (\tilde L)) =b+\lambda_n+\gamma_p(\alpha,c)$. This  proves the upper estimate in (\ref{bestLpJ}). Concerning the lower estimate we recall the bound
$$
\|T_{p,J}(t)\|_p \le M^{\big |1-\frac2p\big |}e^{-(\lambda_n+\gamma_p(\alpha,c))\, t}
$$
proved in Proposition \ref{analyt} for the semigroup generated by $\tilde L$ in $L^p_{J}(\R^N)$ (see  (\ref{expdecay2}) with $\Sigma=S_{N-1}$). Since the resolvent at $b$  is the integral from $0$ to $\infty$ of the generated semigroup, integrating the above bound we get the lower estimate with $c=M^{-1}$.
\qed

\begin{os} {\rm If the condition  $b+\lambda (P_n)+\gamma_p(\alpha,c)>0$ is satisfied in a certain range $[p_1,p_2]$ it is easy to see that the conclusion is still valid in $[p_1,p_2]$. }\end{os}

The description of the spaces $L^p_J(\R^N)$ is simple, when the $L^2(S_{N-1})$ projection onto $F_{J,2}$ extends to a bounded operator in $L^p(S_{N-1}$, see Lemma \ref{projection} (this is always true when $J$ or $\N_0 \setminus J$ is finite). In this case $L^p_J(\R^N)$ is characterized by equation (\ref{caratterizzazione}) as
$$ L^p_{J}({\cal C}_\Sigma)=\left \{u \in L^p({\cal C}_\Sigma): \int_\Sigma u(\rho\, \omega)P_j(\omega) \, d\sigma (\omega)=0\ {\rm for} \ \rho>0\ {\rm  and}\  j \not \in J\right \}.
$$
Special situations of interest are the subspaces $L^p_n$, $L^p_{\ le n}$ and $L^p_{\ge n}$ introduced in Definition \ref{Lp<n}.

According to this notation it follows from Theorem \ref{RellichFJ} that  if $p=1, N/2$ and $\alpha=0$ Rellich inequalities hold in $L^p_{\ge 1}$, that is for every smooth function in $u(\rho, \omega) \in L^1$ such that or every $\rho>0$
$$\int_{S_{N-1}}u(\rho, \omega) \, d\sigma(\omega)=0.
$$  
If $p=N$ and $\alpha=0$, Rellich inequalities hold in $L^p_{\neq 1}$. If $p=\infty, \alpha=0$, we see that Rellich inequalities hold in $L^p_{\neq 2}$. Similar remarks holds for $\alpha \neq 0$, see subsection 7.4.

We can now complement Theorem \ref{Cald-Zyg} showing that weighted Calder\'on-Zygmund inequalities can hold in $L^p_J(\R^N)$ even though they do not hold in $L^p(\R^N)$. We may assume that
$$\alpha \neq 1-\frac{N}{p},\quad \alpha \neq2-\frac{N}{p},\quad \alpha \neq1+\frac{N}{2}-\frac{N}{p}$$ 
otherwise Theorem \ref{Cald-Zyg} applies. Note that the last condition implies that the parabola ${\cal P}_{p,\alpha,0}$ is non degenerate.

\begin{teo} \label{Cald-ZygFJ}
Let $N\geq 3$, $1 < p < \infty$, $\alpha\in\R$, and set $\N_{J}=\{n \in N_0: \ {\rm there\ exists\  } j\in J\ {\rm such\ that\  }\lambda_n=\lambda(P_j)\}$. Assume that
$$
\alpha \neq \frac{N}{p'}+n \qquad \alpha \neq -\frac{N}{p}+2-n \quad {\rm for \ every\ } n \in \N_J.
$$
 Then
there exists a positive constant $C=C(\alpha, N, p)$ such that the weighted Calder\'on-Zygmund inequalities
 \begin{equation*} 
\||x|^\alpha D^2 u\|_p \leq C \||x|^{\alpha} \Delta u\|_p
\end{equation*}
hold for every $u\in C_c^\infty(\R^N \setminus\{0\}) \cap L^p_J(\R^N)$.
\end{teo}
{\sc Proof. } As shown in the proof of Theorem \ref{Cald-Zyg},  Rellich inequalities always imply  Calder\'on-Zygmund inequalities. The  former hold beacuse of Theorem \ref{RellichFJ} and since
the condition $$
\lambda_n+\gamma_p(\alpha,0)=\left (\frac{N}{p}-2+\alpha +n\right )\left (\frac{N}{p'}-\alpha+n\right )\neq 0
$$
for every $n \in \N_J$ is equivalent to that in the statement.
\qed

\begin{os}
{\rm We point out that,  for every $\alpha\in\R$, we can choose  $I$ finite such that  Rellich and Calder\'on-Zygmund inequalities are true in $L^p_{J}$, $J=\N_0\setminus I$.}
\end{os}

Finally, let us consider  Rellich and Calder\'on-Zygmund inequalities for even and odd functions. 
Let $\lambda_n$ be the eigenvalues of $-\Delta_0$ and let $E_n=Ker (\lambda_n+\Delta_0)$, $n \in \N_0$. Then $E_0$ consists of constant functions and each $E_n$ with $n \ge 1$ has dimension greater than 1. 
Let $P_j$, $Q_j$  be odd and even spherical harmonics, respectively.
Then the set $\{\lambda(P_j)\}$ coincides with $\{\lambda_n, n \in \N\}$ and  $\{\lambda(Q_j)\}$ coincides with $\{\lambda_n, n \in \N_0\}$. From Theorem \ref{RellichFJ} we deduce

\begin{prop} \label {evenodd}
Rellich inequalities
\begin{equation*} 
\||x|^\alpha Lu\|_p \geq C\||x|^{\alpha-2} |u|\|_p
\end{equation*}
hold for smooth even functions if and only if $b+\lambda_n \not \in {\cal P}_{p,\alpha,c}$ for every $n\in N_0$ and for odd functions if and only if $b+\lambda_n \not \in {\cal P}_{p,\alpha,c}$ for every $n\in N$ .
\end{prop}

Observe that odd and even functions constitute complemented subspace of $L^p(\R^N)$ and that the $L^2$ projections on even and odd functions on the sphere extend to bounded operator on $L^p(S_{N-1})$, so that characterization (\ref{caratterizzazione}) hold.

Calder\'on-Zygmund inequalities for odd functions follow from Rellich inequalities, as above, see also the next subsection.

\subsection{Rellich and Calder\'on-Zygmund inequalities in $\R^N_+$}
By using the results of Section 6 on the spectrum of $A_p^+$  in $L^p(\R^N_+)$ 
we  can improve Rellich and Calder\'on-Zygmund inequalities in $\R^N_+$ for smooth functions vanishing at the boundary. The following result is proved exactly as Theorem \ref{Rellich-p} and differs from it since the condition $b+\lambda_n \not \in {\cal P}_{p,\alpha,c}$ is required for $n \in \N$ and not for $n \in \N_0$.

\begin{teo} \label{Rellichhalf}
Let $N\geq 2$, $\alpha, b,\ c\in\R$, $1 \le p \le \infty$. Let $\gamma_p(\alpha,c)$ and ${\cal P}_{p,\alpha,c}$ as in (\ref{gammap}) and (\ref{Pp}).
Then there exists a positive constant $C=C(N,\alpha,  p, c, b)$ such that
\begin{equation*} 
\||x|^\alpha Lu\|_p \geq C\||x|^{\alpha-2} |u|\|_p
\end{equation*}
holds for every  for every $u\in C_c^\infty(\overline{\R^N_+}\setminus\{0\})$ 
satisfying $u=0$ on $\partial \R^N_+$  if and only if $b+\lambda_n \not \in {\cal P}_{p,\alpha,c}$ for every $n \in \N$.
\end{teo}
{\sc Proof.} Proceed as in Theorem \ref{Rellich-p} considering the operator
$$ 
\tilde L=|x|^2\Delta +(4-2\alpha+c)x\cdot\nabla +(2-\alpha)(N-\alpha+c).
$$
in the half-space, with Dirichlet boundary conditions, whose spectrum is computed in Proposition \ref{specRN+}.
\qed

\begin{prop} \label{besthalf}
If $b+\gamma_p(\alpha,c)+\ds\frac{4(p-1)(N-1)}{p^2}>0$, then the constant $C$ satisfies 
$$b+\gamma_p(\alpha,c)+ \ds\frac{4(p-1)(N-1)}{p^2}\leq C \leq b+\gamma_p(\alpha,c)+\lambda_1=b+\gamma_p(\alpha,c)+N-1.$$ 
\end{prop}
{\sc Proof. }Indeed the lower bound of $C$ follows from the dissipativity estimate of $\tilde{L}$ given in Proposition \ref{dissipativity2} as in the proof of Theorem \ref{Rellich-p} and the upper bound from the computation of the spectral bound, see Proposition \ref{specRN+} since, as in Theorem \ref{Rellich-p}, $C \le dist (b, \sigma(\tilde L))$.
\qed
Observe that the upper and the lower bound coincide, in the above proposition, if and only if $p=2$. We do not know the exact value of $C$, however we can prove its asymptotic behavior  as $b \to \infty$.

\begin{prop} \label{nobesthalf}
We have 
$$C -\left ( b+\gamma_p(\alpha,c)+\ds\frac{4(p-1)(N-1)}{p^2}\right ) \to 0$$
as $b \to +\infty$
\end{prop}
{\sc Proof.} We may assume that $p \neq 2$. Recall that, from the proof of Theorem \ref{Rellich-p}, $C=\|(b-{\tilde L})^{-1}\|^{-1}_p$, 
where $\tilde L$ is considered in the half-space with Dirichlet boundary conditions. Let 
$$
t_p=-\gamma_p(\alpha,c)-\ds\frac{4(p-1)(N-1)}{p^2} \ge s_p=-\gamma_p(\alpha,c)-(N-1)
$$
($s_p$ is the spectral bound of $\tilde L$) and, given $\eps>0$, assume that the best constant $C_b$ satisfies $C_{b_n} \ge b_n-t_p+\eps$ for a sequence $b_n \to \infty$.
Hence $\|(b_n-{\tilde L})^{-1}\|_p\le(b_n-t_p+\eps)^{-1}$ for every $n$ and, by the resolvent equation, 
 $\|(b-{\tilde L})^{-1}\|_p \le (b-t_p+\eps)^{-1}$ for every $b >t_p-\eps$.
This implies the estimate $\|T^+(t)\|_p \le e^{(t_p -\eps) t}$ for the generated semigroup, which however, contradicts Proposition \ref{dissipativity2}, see also Remark \ref{noncontractivity}. Then $C_b < b-t_p +\eps$ for large $b$ and, since $b-t_p \le C_b$, by Proposition \ref{besthalf}, the proof is complete.
\qed

Finally, let us consider  Calder\'on-Zygmund inequalities in $\R^N_+$ for smooth functions vanishing at the boundary. If Calder\'on-Zygmund inequalities hold in $\R^N$ for a certain $\alpha \in \R$, then they hold in the half space for the same $\alpha$ by a simple reflection argument (for a function $u$  vanishing at the boundary, consider its odd reflection with respect the axis $x_N$). By Theorem \ref{Cald-Zyg} it follows that the above inequalities are satisfied in the half-space if $N \ge 3$ and
$$
\alpha \neq \frac{N}{p'}+n \quad {\rm for\  every\ } n \ge 0 \quad {\rm and\  }\quad \alpha \neq -\frac{N}{p}+2-n \quad {\rm for\  every\ } n \ge2.
$$
 for every $n \in \N_0$. However Calder\'on-Zygmund inequalities fail in the whole space for $\alpha=N/p'$ (corresponding to $\lambda_0=0$ and to radial functions). Let us show that this value is permitted in the half-space and also for $N=2$.

\begin{teo} \label{CZhalf}
Let $N\geq 2$, $1 < p < \infty$, $\alpha=\frac{N}{p'}$. Then
there exists a positive constant $C=C(\alpha, N, p)$ such that the weighted Calder\'on-Zygmund inequalities
 \begin{equation*} 
\||x|^\alpha D^2 u\|_p \leq C \||x|^{\alpha} \Delta u\|_p
\end{equation*}
hold for every  for every $u\in C_c^\infty(\overline{\R^N_+}\setminus\{0\})$ 
satisfying $u=0$ on $\partial \R^N_+$.
\end{teo}
{\sc Proof.} It is sufficient to show that Rellich inequalities hold and this follows from Theorem \ref{Rellichhalf} setting $b=c=0$. The condition $\lambda_n \not \in {\cal P}_{p,\alpha,0}$ for $n \in \N$ reads
$$
\lambda_n +\left (\frac{N}{p}-2+\alpha \right )\left (\frac{N}{p'}-\alpha \right ) \neq 0
$$
for every $n \in \N$ and is clearly satisfied since $n=0$ is omitted.
\qed

\begin{os} \label{odd}
{\rm Observe that the same phenomenon occurs for odd functions of the preceeding subsection. In fact, since $W^{2,p}(\R^N_+)\cap W^{1,p}_0(\R^N_+)$ can be identified with the subspace of $W^{2,p}(\R^N)$ consisting of odd functions, all the results in this subsection can be formulated for odd functions in $\R^N$ rather for functions on $\R^N_+$ vanishing at the boundary. Observe that this identification is obtained via odd reflection of a function $u$  with respect to the $x_N$ axis which has the effect of multiplying by $2^{1/p}$ the norm of $u, \nabla u, D^2 u$.}\end{os}

\subsection{Best constants on special subspaces}
Here we specialize Theorem \ref{RellichFJ} to the case of subspaces generated by spherical harmonics of a fixed order, computing also the best constants. For simplicity we assume $\Sigma=S_{N-1}$ so that $L^p({\cal C}_\Sigma)=L^p(\R^N)$. For $n \in \N_0$ let, according to Definition \ref{Lp<n}, $L^p_n=\{u=\sum_{i=1}^{d_n}f_i(\rho)P_i(\omega)\, \}$,
where $f_i \in L^p((0,\infty), \rho^{N-1}\, d\rho)$ and  $(P_i)$ is a basis for the space of spherical harmonics of order $n$, whose dimension is $d_n$. Note that $d_0=1$ and that $L^p_0$ consists of radial functions. Note also that $L^p_n$ is closed, by Lemma \ref{projection}.

\begin{teo} \label{bestconstant} 
Let $N\geq 2$, $\alpha, b,\ c\in\R$, $1 \le p \le \infty$.
Then there exists a positive constant $C=C(N,\alpha,  p, c, b)$ such that
$$
\||x|^\alpha Lu\|_p \geq C\||x|^{\alpha-2} |u|\|_p
$$
holds for every $u\in C_c^\infty(\R^N \setminus \{0\})\cap L^p_n$  if and only if $b+\lambda_n \not \in {\cal P}_{p,\alpha,c}$.
   If, in addition,  $b+\lambda_n+\Bigl(\frac{N}{2}-1+\frac{c}{2}\Bigr)^2 \ge0$
the optimal constant is given by
$C=|b+\lambda_n+\gamma_p(\alpha,c)|.$
\end{teo}
{\sc Proof.} The first statement follows from Theorem \ref{RellichFJ}, if we select $J$ corresponding to all spherical harmonics of order $n$. In order to prove the second part, however, we repeat some arguments of the proof.
As in Theorem \ref{Rellich-p} we set $v(x)=|x|^{\alpha-2}u(x)=\sum_{i=1}^{d_n}\rho^{\alpha-2}f_i(\rho)P_i(\omega)$ and observe that
$$|x|^\alpha L u:=\tilde{L}v-(b+\lambda_n)v$$
where $\tilde L=\Gamma \otimes I$ and $\Gamma$ is the radial operator 
$$
\Gamma =\rho^2 D_{\rho \rho} +(N+3-2\alpha+c)\rho D_\rho +(2-\alpha)(N-\alpha+c).
$$
Rellich inequality in the statement is equivalent to
$$
\|\tilde{L}v-(b+\lambda_n)v\|_p \geq C\|v\|_p,
$$
with the same constant $C$.
Next observe that ${\tilde L}-(b+\lambda_n)=\left (\Gamma -(b+\lambda_n )\right)\otimes I$ so that $b+\lambda_n$ belongs to the spectrum of $\tilde L$ if and only if it belongs to the spectrum of $\Gamma$
which is the parabola (\ref{Pp}), by Proposition \ref{spec-rad}.

Moreover 
$$
C=\|\left({\tilde L}-(b+\lambda_n)\right)^{-1}\|_p=\|\left(\Gamma-(b+\lambda_n)\right)^{-1}\|_p=|b+\lambda_n+\gamma_p(\alpha,c)|
$$
if  $b+\lambda_n+\Bigl(\frac{N}{2}-1+\frac{c}{2}\Bigr)^2>0$, by Propositions \ref{normagamma}, \ref{tensornorm} (see also Remark \ref{tensornorm1}).
\qed

Let us specialize the above result to the case of the Laplacian. Note that the case of radial functions, corresponding to $n=0$, is already contained in \cite[Theorem 1.1]{musina1}.

\begin{cor} \label{bc1}
Assume that $b=c=0$, that is $L=\Delta$. 
Then Rellich inequalities hold for smooth functions in $L^p_n$  if and only if $\alpha \not =N/p'+n$, $\alpha \not =-N/p+2-n$. The optimal constant is given by
$$
C_n(N,\alpha,p)=\Bigl|\Bigl (\frac{N}{p}-2+\alpha+n\Bigr)\Bigl(\frac{N}{p'}-\alpha+n\Bigr)\Bigr|.
$$
\end{cor}
{\sc Proof. } 
The parabola ${\cal P}_{p,\alpha,0}$ degenerates if and only if $\bar{\alpha}=N(1/2-1/p)+1$ and $\gamma_p(\bar{\alpha},0)>0$ if and only if $N>2$. However, if $N=2$, then $\bar {\alpha}=2/p'$, $\gamma_p(\bar{\alpha},0)=0$.  Hence Rellich inequalities hold for $\bar {\alpha}$ for every $n \in \N_0$  if  $N\ge 3$ and for every $n \in \N$  if  $N=2$.
Assume now that $\alpha \neq \bar{\alpha}$.
Since 
$$
\lambda_n+\gamma_p(\alpha,0)=\left (\frac{N}{p}-2+\alpha+n\right )\left (\frac{N}{p'}-\alpha+n\right ),
$$
the condition $\lambda_n+\gamma_p(\alpha,0) \neq 0$  is equivalent to $\alpha \neq N/p'+n$, $\alpha \neq  -N/p+2-n$, as in the statement. 
Finally, since the condition $\lambda_n +(N/2-1)^2 \ge 0$ is clearly satisfied, Theorem \ref{bestconstant} yields $C=|\lambda_n+\gamma_p(\alpha,0)|$, as in the satement.
\qed

When $\alpha=b=c=0$ we know that the classical Rellich inequalities hold except for $p=1,N/2$ where they fail on radial functions, for $p=N$ where they fail on $L^p_1$ and for $p=\infty$ where they fail on $L^p_2$. The previous corollary yields the best best constants on $L^p_0$, $L^p_1$ and $L^p_2$.

\begin{cor} \label{cor2} If $\alpha=b=c=0$ the best constants for Rellich inequalities in $L^p_0$, $L^p_1$ and $L^p_2$ are given by
$$
C_0=\Bigl|\Bigl (\frac{N}{p}-2\Bigr)\Bigl(\frac{N}{p'}\Bigr)\Bigr|,  \quad  C_1=N^2\left (\frac{1}{p'} +\frac1N \right )\left| \frac1p-\frac1N\right |, \quad  C_2=\frac{N}{p} \left (\frac{N}{p'}+2 \right ).
$$
\end{cor}

Note that $C_0=0$ when $p=1,N/2$, $C_1=0$ when $p=N$, $C_2=0$ when $p=\infty$.
\begin{os} \label{sphericalp}
{\rm It is important to note that Theorem \ref{Rellich-p} can be deduced from Theorem \ref{bestconstant} only if $p=2$ as done in Section 2. This happens because spherical harmonics are an orthogonal basis of $L^2(S_{N-1})$ but are not even a Schauder basis of $L^p(S_{N-1})$ when $N \ge 3$. In particular, keeping the notation of Theorem \ref{bestconstant}, the inequality
\begin{equation} \label{noeq}
C(N,\alpha,p,b,c) \leq\min_{n \in \N_0}|b+\lambda_n+\gamma_p(\alpha,c)|
\end{equation} for the best constant in Rellich inequalities holds but equality can fail for $p \neq 2$.
An example of this phenomenon is exhibited in Proposition \ref{nobesthalf}, though formulated in the half-space. Note that $L^p(\R^N_+)$ can be regarded as a subspace of $L^p_{>0}$  by extending any function vanishing at the boundary to an odd function in $\R^N$. If equality were true in (\ref{noeq}), by applying it to $L^p_{>0}$ and large positive $b$ we would obtain $C=b+\lambda_1+\gamma_p(\alpha,c)$ also in the half-space, contradicting Proposition \ref{nobesthalf}.
}
\end{os}
\subsection{Remarks on best constants} \label{estimatebest}
In this section we assume that $L=\Delta$, that is $b=c=0$. As already pointed out in the introduction, the best constants for Rellich inequalities (\ref{WRI}) are known only for $N\geq 3$ and for $2-\frac{N}{p}<\alpha<\frac{N}{p'}$. They are given by $(N/p-2+\alpha)(N/p'-\alpha) $ and coincide with the best constants on radial funcions, see the previous subsection. The best constants for $\alpha$ outside the above range seem to be unknown even for $\alpha=0$. Note however, that if $\alpha=0$ and $p>N/2$ the best constant on the whole space can be strictly larger than the corresponding one on radial functions since Rellich inequalities fail for $p=N$ (on $L^p_1$) but hold on radial functions. Similar examples can be done for every $\alpha$.

The following result yields an estimate of the $p$-dependence of the best constant on $L^p_{\ge n}$. Note, however, that the constant $c$ below  is not explicitely given.

\begin{cor} \label{bestgen}
Let $N \ge 2$, $\alpha \in R$, $b=c=0$ and $n \ge 1$. If $\left (\frac{N}{p}-2+\alpha+n\right )\left (\frac{N}{p'}-\alpha+n\right )>0$ for  every $1 \le p \le \infty$, then the best constants $C_{\ge n}(N,\alpha,p)$ for Rellich inequalities in $L^p_{\ge n}$ satisfy
$$
c^{\big |1-\frac2p |}\left (\frac{N}{p}-2+\alpha+n\right )\left (\frac{N}{p'}-\alpha+n\right ) \le C_{\ge n}(N,\alpha,p) \le 
\left (\frac{N}{p}-2+\alpha+n\right )\left (\frac{N}{p'}-\alpha+n\right ),
$$
where $c>0$ depends on $n,N$ but not on $p$.
\end{cor}
{\sc Proof. } This follows immediately from (\ref{bestLpJ}) of Theorem \ref{RellichFJ} since $b,c=0$ and $$\lambda_n+\gamma_p(\alpha,0)=\left (\frac{N}{p}-2+\alpha+n\right )\left (\frac{N}{p'}-\alpha+n\right ).
$$
\qed

An estimate of the best constants on the whole space $L^p(\R^N)$ can be given by combining Corollaries \ref{bc1},   \ref{bestgen}.

\begin{prop}  \label{bestgenp} Let $N \ge 2$, $\alpha \in R$, $b=c=0$ and $C(N,p,\alpha)$ be the best constants for Rellich inequalities (with the understanding that Rellich inequalities do not hold if $C(N,p,\alpha)=0$). Then
$$
c\min_{n \in \N_0}\left |\left (\frac{N}{p}-2+\alpha+n\right )\left (\frac{N}{p'}-\alpha+n\right ) \right |\le C(N, \alpha,p) \le 
\min_{n \in \N_0}\left |\left (\frac{N}{p}-2+\alpha+n\right )\left (\frac{N}{p'}-\alpha+n\right )\right |,
$$
where $c$ depends on $\alpha, N$ but not on $1 \le p \le \infty$.
\end{prop}

{\sc Proof. }  The upper estimate has been already observed in Remark \ref{sphericalp}, see (\ref{noeq}). To prove the lower estimate, let $n \in \N_0$ be the first integer such that $\lambda_k +\gamma_p(\alpha,0)>0$ for every $k \ge n$ and every $1\le p\le \infty$.
Let $P_i$, $1=1, \dots, n-1$, be the orthogonal projections from $L^p(\R^N)$ onto $L^p_i$, which are bounded in $L^p(\R^N)$ uniformly for $1 \le p \le \infty$. If $Q=P_0+\cdots P_{n-1}$ and $P=I-Q$, then $Q$ and $P$ are the orthogonal projections onto $L^p_{<n}$ and $L^p_{\ge n}$, respectively. Then 
$$
L^p(\R^N)=L^p_0\oplus \cdots \oplus L^p_{n-1} \oplus L^p_{\ge n}
$$
 Let $\kappa>0$ be such that $\|P_i\|_p \le \kappa$ for every $i=1, \dots, n$, $1 \le p\le \infty$.  If $u \in L^p(\R^N)$, then $u=\sum_{i=0}^n P_i u$.  From Corollaries \ref{bc1}, \ref{bestgen} we have for $u \in C_c^\infty (\R^N \setminus \{0\})$
$$
C_i \||x|^{\alpha-2}P_i u\|_p \le \||x|^\alpha \Delta P_i u\|_p =\|P_i(|x|^\alpha \Delta  u)\|_p \le \kappa \||x|^\alpha \Delta  u\|_p
$$
where $C_i=c_i\left |\left (\frac{N}{p}-2+\alpha+i\right )\left (\frac{N}{p'}-\alpha+i\right )\right |$ and $c_i=1$ for $i<n$, $c_n=c$.
If $C=\min_{i=1, \dots n}{C_i}$ then 
$$
C \||x|^{\alpha-2} u\|_p \le C\sum_{i=1}^n  \|P_i(|x|^{\alpha-2} u)\|_p =C\sum_{i=1}^n  \||x|^{\alpha-2}P_i u\|_p \le  n \kappa \||x|^\alpha \Delta  u\|_p
$$
and the thesis follows with $c=C/(n \kappa)$.
\qed

In particular, for $\alpha=0$ we obtain
\begin{equation} \label{best0}
\left\{
\begin{array}{ll} \displaystyle
C(N,p)= \left (\frac{N}{p}-2\right )\left (\frac{N}{p'}\right  ) \ {\rm if \ }\ 1<p<\frac{N}{2}\\ 
\displaystyle  c\min_{n =0,1} \left |\left (\frac{N}{p}-2+n\right )\left (\frac{N}{p'}+n\right ) \right |\le C(N,p) \le 
\min_{n=0,1}\left |\left (\frac{N}{p}-2+n\right )\left (\frac{N}{p'}+n\right )\right | \ {\rm if \ }\ \frac{N}{2} <p<N\\
 \displaystyle
c\min_{n =1,2} \left |\left (\frac{N}{p}-2+n\right )\left (\frac{N}{p'}+n\right ) \right |\le C(N,p) \le 
\min_{n=1,2}\left |\left (\frac{N}{p}-2+n\right )\left (\frac{N}{p'}+n\right )\right | \ {\rm if \ }\ N<p<\infty
\end{array}\right.
\end{equation}

for the best constant $C(N,p)$ of the unweighted Rellich inequalities. A more explicit estimate of the lower bound above can be obtained by exploiting the relationships between Rellich and Calder\'on-Zygmund inequalities, as in Section 4. Setting $\alpha=0$ in equation (\ref{hardytwice}) we obtain
$$
\left |\left (\frac{N}{p}-2\right )\left (\frac{N}{p}-1 \right )\right | \||x|^{-2} u\|_p \leq  \| D^2 u\|_p \le C_Z(N,p)\|\Delta u\|_p
$$
where $C_Z(N,p)$ is the best constant in  Calder\'on-Zygmund inequalities. Since $C_Z(N,p) \le c_N \frac{p}{p-1}$, where $c_N$ can be estimated through the constants in the weak $L^1$-estimate of the Calder\'on-Zygmund kernel and those of the Marcinkiewicz interpolation  theorem, see \cite[6.2, Chapter II]{stein-weiss},  a lower estimate, similar to that of  (\ref{best0}) but perhaps more explicit follows. Observe also that this argument can be reversed: having a bound on the constants in Rellich inequalities one can obtain, for $p\neq N/2,N$ a bound for $C_Z$ through the constants in Hardy inequalities and the constant $C$ in Lemma \ref{interpolative}, see Section 4.

\section*{Appendix A: Weighted Hardy inequalities in $L^p(\R^N)$}

In this appendix we state and proof three inequalities of Hardy type we need in the paper. The first two are well-known in the literature but we give a short proof for completeness. Proposition \ref{hardy2} (ii) seems to be new, concerning the computation of the best constant.

\begin{prop}  \label{hardy}
Let $1<p<\infty$ and $\beta\in\R$. 
Then $\beta-p+N\neq 0$ if and only if $u\in C_c^\infty(\R^N\setminus\{0\})$, 
\begin{equation}\label{WHardy}
\int_{\R^N}
|x|^{\beta}|\nabla u|^p
\,dx
\geq 
C\int_{\R^N}
|x|^{\beta-p}|u|^p
\,dx
\end{equation}
holds with the optimal constant 
\[
C=\left|
  \frac{\beta-p+N}{p}
\right|^p>0.
\]
\end{prop}

{\sc Proof.}
Let $u\in C_c^\infty(\R^N\setminus\{0\})$ and $\beta-p+N\neq 0$. 
Then integration by parts and H\"older inequality imply
\begin{align*}
\int_{\R^N}|x|^{\beta-p}|u|^p\,dx
&=\frac{-p}{\beta-p+N}\int_{\R^N}|x|^{\beta-p}x\cdot |u|^{p-2}{\rm Re}(\overline{u}\nabla u)\,dx\\
&\leq\frac{p}{|\beta-p+N|}\int_{\R^N}|x|^{\beta-p+1}|u|^{p-1}|\nabla u|\,dx\\
&\leq\frac{p}{|\beta-p+N|}\Bigl(\int_{\R^N}|x|^{\beta-p}|u|^{p}\,dx\Bigr)^{1-\frac{1}{p}}
\Bigl(\int_{\R^N}|x|^{\beta}|\nabla u|^p\,dx\Bigr)^\frac{1}{p}.
\end{align*}
Hence we have \eqref{WHardy}.

Next we show the optimality. 
Fix $\phi\in C_c^\infty(]0,\infty[)$ such that $0\leq\phi\leq 1$ and $\phi(\rho)=1$ on $[\frac{1}{2},2]$.
Choosing the sequence $\{u_{\varepsilon,m}\}_{\varepsilon,m}\subset C_c^\infty(\R^N\setminus\{0\})$ as
\[
u_{\varepsilon,m}(x)=\Gamma(p\varepsilon)^{-\frac{1}{p}}\phi(|x|^{\frac{1}{m}})|x|^{-\frac{\beta-p+N}{p}+\varepsilon}e^{-\frac{|x|}{p}},
\]
we have
\begin{align*}
\int_{\R^N}|x|^{\beta-p}|u_{\varepsilon,m}|^p\,dx
&\to
1\quad(m\to \infty), 
\\
\int_{\R^N}|x|^{\beta}|\nabla u_{\varepsilon,m}|^p\,dx
&\to
\Gamma(p\varepsilon)^{-1}\int_0^\infty
\left|\frac{\beta-p+N}{p}-\varepsilon+\frac{1}{p}\rho\right|^p\rho^{p\varepsilon-1}e^{-\rho}\,d\rho
\quad(m\to \infty).
\end{align*}
Noting that \eqref{WHardy} and 
\[
\Bigl(
\Gamma(p\varepsilon)^{-1}
\int_0^\infty
\left|\frac{\beta-p+N}{p}-\varepsilon+\frac{1}{p}\rho\right|^p\rho^{p\varepsilon-1}e^{-\rho}\,d\rho
\Bigr)^{\frac{1}{p}}
\leq
\left|\frac{\beta-p+N}{p}-\varepsilon\right|
+\frac{1}{p}\Bigl(\frac{\Gamma(p+p\varepsilon)}{\Gamma(p\varepsilon)}\Bigr)^{\frac{1}{p}},
\]
we see that
\[
\lim_{\varepsilon\downarrow0}\left(\lim_{m\to\infty}
\int_{\R^N}|x|^{\beta}|\nabla u_{\varepsilon,m}|^p\,dx
\right)
=
\left|\frac{\beta-p+N}{p}\right|^p,
\]
thus showing the optimality of \eqref{WHardy}.
\qed

\begin{os} {\rm If $\beta-p+N=0$, then $u_{\varepsilon,m}$ is nothing but a counterexample to 
\eqref{WHardy}.}
\end{os}

\begin{prop}\label{hardy2}
Let $1<p<\infty$ and $\beta\in\R$. Then the following inequalities hold:
\medskip

\noindent
(i) if $N-2+\beta\neq 0$, then
for every $u\in C_c^\infty(\R^N\setminus\{0\})$, 
\begin{equation}\label{WHardy2-RN}
\int_{\R^N}
|x|^{\beta}|\nabla u|^{2}|u|^{p-2}
\,dx
\geq 
\left(\frac{N-2+\beta}{p}\right)^2
\int_{\R^N}
|x|^{\beta-2}|u|^{p}
\,dx;
\end{equation}
(ii) for every $u\in C_c^\infty(\overline{\R^N_+}\setminus\{0\})$ 
satisfying $u=0$ on $\partial \R^N_+$,
\begin{equation}\label{WHardy2-half}
\int_{\R^N_+}
|x|^{\beta}|\nabla u|^{2}|u|^{p-2}
\,dx
\geq 
\left[\left(\frac{N-2+\beta}{p}\right)^2+\frac{4(N-1)}{p^2}\right]
\int_{\R^N_+}
|x|^{\beta-2}|u|^{p}
\,dx.
\end{equation}
The constants in \eqref{WHardy2-RN} and \eqref{WHardy2-half} are sharp.
\end{prop}

{\sc Proof.}
Let $\Omega$ be $\R^N\setminus\{0\}$ or $\R^N_+$. 
First we show \eqref{WHardy2-RN} and \eqref{WHardy2-half} for $u\in C_c^\infty(\Omega)$. Set a function 
$Q\in C^2(\Omega)$ such that $Q>0$ 
($Q$ will be chosen depending on $\Omega=\R^N$ or $\Omega=\R^N_+$). 
Put $v(x):=Q(x)^{-\frac{2}{p}}u(x)$. Then noting that 
\[
|\nabla u|^2=Q^{\frac{-2(p-2)}{p}}
\left(|\nabla v|^2+\frac{4}{p}Q\nabla Q{\rm Re}(\overline{v}\nabla v)+\frac{4}{p^2}|\nabla Q|^2|v|^2\right)
\]
we see from integration by parts that 
\begin{align*}
\int_{\Omega}
|x|^{\beta}|\nabla u|^{2}|u|^{p-2}
\,dx
&\geq 
\frac{4}{p}
\int_{\Omega}
|x|^{\beta}Q\nabla Q\cdot{\rm Re}(\overline{v}\nabla v)|v|^{p-2}\,dx
+
\frac{4}{p^2}
\int_{\Omega}
|x|^{\beta}|\nabla Q|^2|v|^p\,dx
\\
&=
-
\frac{4}{p^2}
\int_{\Omega}
|x|^{\beta}\left(\Delta Q+\frac{\beta x\cdot\nabla Q}{|x|^2}\right)Q|v|^p\,dx.
\\
&= 
-
\frac{4}{p^2}
\int_{\Omega}
|x|^{\beta}\left(\Delta Q+\frac{\beta x\cdot\nabla Q}{|x|^2}\right)Q^{-1}|u|^p\,dx.
\end{align*}
{(Case (i): $\Omega=\R^N\setminus\{0\}$)} We choose $Q_{0}(x)=|x|^{1-\frac{N+\beta}{2}}$. Then by easy culculation we have 
\[
\Delta Q_{0} + \frac{\beta x\cdot\nabla Q_{0}}{|x|^2}
=-\left(\frac{N-2+\beta}{2}\right)^2\frac{Q_{0}}{|x|^2}.
\]
Thus we obtain 
\begin{equation}
\int_{\R^N}
|x|^{\beta}|\nabla u|^{2}|u|^{p-2}
\,dx
\geq 
\frac{4}{p^2}
\left(\frac{N-2+\beta}{2}\right)^2
\int_{\R^N}
|x|^{\beta-2}|u|^{p}\,dx.
\end{equation}
We have proved (i).
\\
{(Case (ii): $\Omega=\R^N_+$)} We choose $Q_{+}(x)=x_n|x|^{-\frac{N+\beta}{2}}$. Then similary we have 
\[
\Delta Q_{+} + \frac{\beta x\cdot\nabla Q_{+}}{|x|^2}
=-\left[\left(\frac{N-2+\beta}{2}\right)^2+N-1\right]\frac{Q_{+}}{|x|^2}.
\]
Thus we obtain 
\begin{equation}
\int_{\R^N}
|x|^{\beta}|\nabla u|^{2}|u|^{p-2}
\,dx
\geq 
\frac{4}{p^2}
\left[\left(\frac{N-2+\beta}{2}\right)^2+N-1\right]
\int_{\R^N}
|x|^{\beta-2}|u|^{p}\,dx. 
\end{equation}
We have proved (ii) for $u\in C_c^\infty(\R^N_+)$. 

Now we prove 
(ii) for $u\in C_c^\infty(\overline{\R^N_+}\setminus\{0\})$ satisfying $u=0$ on $\partial \R^N_+$. We introduce $G\in C^\infty(\R)$ 
and $G_n\in C^\infty(\R)$ for $n\in \N$ as 
\[
G(s)=
\begin{cases}
1 & {\rm if\ } |s|\geq 2,
\\
0 & {\rm if\ } |s|<1
\end{cases}
\]
and $G_n(s)=G(ns)$. Then applying \eqref{WHardy2-half} to $u_n(x)=G_n(|u|)u\in C_c^\infty(\R^N_+)$,
we have 
\[
\int_{\R^N_+}
|x|^{\beta}|\nabla u_n|^{2}|u_n|^{p-2}
\,dx
\geq 
\left[\left(\frac{N-2+\beta}{p}\right)^2+\frac{4(N-1)}{p^2}\right]
\int_{\R^N_+}
|x|^{\beta-2}|u_n|^{p}
\,dx.
\]
Noting that $\supp u_n\subset \supp u$, 
$u_n\to u$ uniformly on $\supp u$
and $\nabla u_n \to \nabla u$ in $L^2(\supp u)$, 
from the dominated convergence theorem  we obtain \eqref{WHardy2-half}.
\qed

\begin{os}
{\rm The optimality of \eqref{WHardy2-RN} and \eqref{WHardy2-half} also can be given by using the 
following sequences, respectively:
\[
u_{\varepsilon,m}(x)
=
\begin{cases}
\phi_m(|x|)|x|^{\frac{2-N-\beta}{p}+\varepsilon} e^{-\frac{|x|}{p}} 
&{\rm if}\ \Omega=\R^N\setminus\{0\},
\\[3pt]
\phi_m(|x|)
x_N^{\frac{2}{p}}|x|^{-\frac{N+\beta}{p}+\varepsilon} e^{-\frac{|x'|}{p}} 
&{\rm if}\ \Omega=\R^N_+,
\end{cases}
\]
where $|x'|=(x_1^2+\cdots+x_{N-1}^2)^{\frac{1}{2}}$ and 
$\{\phi_m\}$ is a suitable family of cut-off functions.}
\end{os}

\section*{Appendix B: The operator $A$ in continuous function spaces}

Let $\Omega=\R^N \setminus\{0\}$.
We consider the operator $A$
endowed with its  maximal domain in $C_0^0(\R^N)$  
\begin{equation*}
D^0_{max}(A)= \{u \in C_0^0(\R^N)\cap W_{loc}^{2,p}(\Omega) {\rm \ for \ every\ }  p<\infty: Au \in C_0^0(\R^N)\}.
\end{equation*}                                                                                                              
We start by  studying existence and uniqueness in the larger space
\begin{equation*}
D_{max}(A)= \{u \in C_b(\Omega)\cap W_{loc}^{2,p}(\Omega) {\rm \ for \ every\ }  p<\infty: Au \in C_b(\Omega)\}.
\end{equation*}  
for the elliptic equation
\begin{equation}  \label{elliptic}
\lambda u-Au=f
\end{equation}
for $\lambda >0$ and  $f\in C_b(\Omega)$.

\begin{prop}  \label{esistenza}
For every $f\in C_b(\Omega)$, $\lambda>0$, there exists $u\in
D_{max}(A)$ solving equation (\ref{elliptic}) and satisfying the
inequality $\|u\|_\infty\leq\|f\|_{\infty}/\lambda$. Moreover, $u
\ge 0$ whenever $f \ge 0$.
\end{prop}
{\sc Proof.} The proof is identical to that given in \cite[Theorem
3.4]{met-pal-wack}. In fact, for every $f\in C_b(\Omega)$,  a solution
$u$ of (\ref{elliptic}) can be obtained as limit of solutions $u_n$ of the
Dirichlet problems associated with the operator above in the
sequence of annuli $B_n\setminus B_{\frac{1}{n}}$ which fill
the whole $\Omega$. 
\qed

Uniqueness  follows from the existence of suitable Lyapunov
functions for the operator $A$ (see \cite[Theorem 3.7]{met-pal-wack}).
\begin{defi}\label{LF}
We say that $V$ is a Lyapunov function for $L$ if $V\in
C^2(\Omega)$, $V\geq 1$, $V\to \infty$ as $|x|\to 0, \infty$ and
$\lambda_0 V-AV\geq 0$ for some $\lambda_0>0$.
\end{defi}

\begin{prop}    \label{uniqueness}
Suppose that there exists $V$  Lyapunov function for the operator
$L$. Then
 $\lambda-A$ is injective on $D_{max}(A)$ for every $\lambda>0$.
\end{prop}

\begin{os} \label{iniettivita}
{\rm Let $0\leq \phi\leq 1$ be a smooth function such that
$\phi(x)=-1$ for $|x|\leq 1/2$ and $\phi(x)=1$ for $|x|\geq 1$. By
easy computations it follows that the function
$V(x)=\phi(x)\ln|x|+1$ is a Lyapunov function for $A$. Therefore
$\lambda-A$ is injective on $D_{max}(A)$ for every $\lambda>0$.}
\end{os}

\begin{prop}
For every $\lambda>0$, the resolvent operator $R(\lambda,A)$ preserves $C_0^0(\R^N)$ that is $\lambda-A$ is bijective from $D^0_{max}(A)$ onto $C_0^0(\R^N)$.
\end{prop}
{\sc Proof.} The proof is as in \cite[Theorem 3.17]{met-pal-wack}) once one shows the existence of a function  $V\in C^2(\Omega)$ such that $V(x)\to 0$ as $|x|\to\infty$ and as $|x|\to 0$, $\lambda_0 V-AV\geq 0$ for some $\lambda_0>0$. For example, we can choose  $V(x)=|x|^2$ in $B_1$, $V(x)=(\log|x|)^{-1}$ in $\R^N\setminus B_2$. It follows that, for $\lambda_0>0$ large enough, $\lambda_0 V-AV\geq 0$.
\qed

By the Hille-Yosida theorem we obtain the following result.

\begin{prop} 
The operator $(A, D^0_{max}(A))$ generates a strongly continuos semigroup of positive contractions in $C_0^0(\R^N)$.
\end{prop}

Now we briefly study the radial operators associated to $A$ and classify the endpoints, according to Feller's theory, see \cite[Section VI.4.c]{engel-nagel}).

\begin{prop} \label{feller}  
Let $N\geq 2$, 
 $$\Gamma =\rho^2D^2+(N-1+c)\rho D$$ in $]0,\infty[$.
Then the origin and infinity are natural points. 
\end{prop}
{\sc Proof.} 
According with the notation used in \cite[Section VI.4.c]{engel-nagel}, we compute the Wronskian  $W$ and the functions $Q$, $R$
\begin{align*}
& W(\rho)=\exp\left\{-\int_1^\rho \frac{N-1+c}{r}dr\right\}=\frac{1}{\rho^{N-1+c}};\\
& Q(\rho)=\frac{1}{\rho^{3-N-c}}\int_1^\rho\frac{1}{r^{N-1+c}}dr=\left\{
\begin{array}{ll}
\displaystyle\frac{1}{2-N-c}\left(\frac{1}{\rho}-\frac{1}{\rho^{3-N-c}}\right) & c\neq 2-N,\\
\displaystyle\frac{1}{\rho}\log \rho   & c=2-N;
\end{array}\right.\\
&R(\rho)=\frac{1}{\rho^{N-1+c}}\int_1^\rho\frac{1}{r^{3-N-c}}dr=\left\{
\begin{array}{ll}
\displaystyle\frac{1}{N-2+c}\left(\frac{1}{\rho}-\frac{1}{\rho^{N-1+c}}\right) & c\neq 2-N,\\
\displaystyle\frac{1}{\rho}\log \rho   & c=2-N.
\end{array}\right.\\
\end{align*}
It follows that $Q\not\in L^1(1,\infty)$ and $R\not\in L^1(1,\infty)$, therefore  infinity  is a natural point. Similarly for $0$. \qed

\begin{os}
{\rm Since both the endpoints $0$ and $+\infty$ are inaccessible and 
$$\frac{1}{\rho^2W(\rho)}=\frac{\rho^{N-1+c}}{\rho^2}=\rho^{N-3-c}\not\in L^1(0,+\infty),$$ 
by \cite[Prop. 6.2 (ii)]{met-pal-wack} we deduce that does not exist an invariant measure for $A$.}
\end{os}

\section*{Appendix C: The norm of the tensor product}
If $X,Y$ are $L^p$-spaces over $G_1, G_2$ we denote by $X\otimes Y$ the algebraic tensor product of $X,Y$, that is the set of all functions $u(x,y)=\sum_{i=1}^n f_i(x)g_i(y)$ where $f_i \in X, g_i \in Y$ and $x \in G_1, y\in G_2$.
If $T,S$ are linear operators  on $X,Y$ we denote by $T\otimes S$ the operator on $X\otimes Y$ defined by 
$$
T\otimes S \left (\sum_{i=1}^n f_i(x)g_i(y)\right)=\sum_{i=1}^n T f_i(x)Sg_i(y).
$$
The following result is probably well-known and we give a proof only for completeness.

\begin{prop} \label{tensornorm}
Let $T:L^p(G_1) \to L^q(G_1)$, $S:L^p(G_2) \to L^q(G_2)$ be bounded operators. Then $T\otimes S$ extends to a bounded operator from $L^p(G_1\times G_2)$ to $L^q(G_1 \times G_2)$ and $\|T\otimes S\|=\|T\|\|S\|$.
\end{prop}
{\sc Proof}
\begin{align*}
&\int_{G_1\times G_2}\Big |\sum_{i=1}^n (Tf_i)(x)(Sb_i)(y)\Big |^q\, dx\, dy=\int_{G_1}dx \int_{ G_2}\Big|\sum_{i=1}^n (Tf_i)(x)(Sb_i)(y)\Big|^q\, dy \\
&=\int_{G_1}dx \int_{ G_2}\Big|S (\sum_{i=1}^n (Tf_i)(x)b_i )\Big|^q\, dy \le \|S\|^q\int_{G_1}dx \int_{ G_2}\Big| \sum_{i=1}^n (Tf_i)(x)b_i \Big|^q\, dy\\
&=\|S\|^q\int_{G_2}dy \int_{ G_1}\Big| \sum_{i=1}^n (Tf_i)(x)b_i \Big|^q\, dx =\|S\|^q\int_{G_2}dy \int_{ G_1}\Big|T( \sum_{i=1}^n f_i b_i(y)) \Big|^q\, dx \\
&\le \|S\|^q\|T\|^q \int_{G_1 \times G_2}\Big |\sum_{i=1}^n f_i(x)g_i(y)\Big|^q\, dx\, dy.
\end{align*}
This gives $\|T\otimes S\| \le \|T\|\|S\|$ (using the density of the above functions in $L^p(G_1\times G_2)$). To show the converse inequality it is sufficient to consider, for $\eps>0$, $f\in L^p(G_1)$ such that $\|f\|_p=1$ and $\|Tf\|_q \ge (1-\eps)\|T\|$, $g\in L^p(G_2)$ such that $\|g\|_p=1$ and $\|Sg\|_q \ge (1-\eps)\|S\|$. Then $f\otimes g$ has norm 1 in $L^p(G_1 \times G_2)$ and $(T\otimes S) (f\otimes g)=(Tf)\otimes (Sg)$ has norm greater than $(1-\eps)^2\|T\|\|S\|$.
\qed

{\begin{os} \label{tensornorm1}
{\rm The above Proposition can be generalized to subspaces of $L^p$-spaces.
Let $E^p(G_i)$, $E^q(G_i)$ be closed subspaces of $L^p(G_i), L^q(G_i)$, $i=1,2$ and let $T:E^p(G_1) \to E^q(G_1)$, $S:E^p(G_2) \to E^q(G_2)$ be bounded operators. Then $T\otimes S$ extends to a bounded operator from (the closure of) $E^p(G_1)\otimes E^p( G_2)$ to (the closure of ) $E^q(G_1) \otimes E^q(G_2)$ and $\|T\otimes S\|=\|T\|\|S\|$. The norm on 
$E^p(G_1)\otimes E^p( G_2)$  and $E^q(G_1) \otimes E^q(G_2)$  is that induced by $L^p(G_1\times G_2)$ and $L^q(G_1 \times G_2)$, respectively. The proof is the same as above.}
\end{os}


\begin{thebibliography}{99}

\bibitem{arendt}
\refer{W. Arendt:} {Gaussian estimates and interpolation of the spectrum in $L^p$,}{Diff. Int. Eq.,}{Vol. 7,n. 5}{(1994), 1153-1168}.

\bibitem{arendt1}
\refer{W. Arendt, K. R\"{a}biger, A. Sourour: }{Spectral properties of the operator equation $AX+XB=Y$, }{Quart. J. Math. Oxford (2), }{45}{ (1994), 133-149}.

\bibitem{caldiroli}
\refer{P. Caldiroli, R. Musina:}{Rellich inequalities with weights,} {Calc. Var. Partial Differential Equations,} {45} {(2012), no. 1-2, 147-164.}

\bibitem{cup-for}
\refer{G. Cupini, S. Fornaro,}{Maximal regularity in $L^p$ for a class of elliptic operators with unbounded coefficients,}{Diff. Int. Eqs.,}{Vol. 17}{(2004), 259-296.}


\bibitem{davies}
\refbook{E. B. Davies: }{Heat Kernels and Spectral Theory, }{
Cambridge University Press, 1989}.

\bibitem{davi-hinz}
\refer{E. B. Davies, A. M. Hinz,}{Explicit constants for Rellich inequalities in $L^p(\Omega)$,}{Math. Z.,}{227 n.3}{(1998), 511-523.}


\bibitem{engel-nagel}
\refbook{K.J. Engel, R. Nagel:} {One parameter semigroups for linear evolutions equations,} {Springer-Verlag, Berlin, (2000)}.


\bibitem{for-lor}
\refer{S. Fornaro, L. Lorenzi:} {Generation results for elliptic operators with unbounded diffusion coefficients in $L^p$ and $C_b$-spaces,} {Discrete and continuous dynamical sistems,} {18} {(2007),  747-772}.

\bibitem{FMPP}
\refer{S. Fornaro, G. Metafune, D. Pallara, J. Pr\"uss:}{$L^p$--theory for some elliptic and
parabolic problems with first order degeneracy at the boundary,}{J. Math.\ Pures Appl.}{87}
{(2007), 367--393.}



\bibitem{kree}
\refer{Paul Kr\'{e}e: }{Sur les multiplicateurs dans $\mathcal{F}L^p$ avec poids, }{Ann. Inst. Fourier, Grenoble, }{Vol. 16, N. 2}{(1966), 91-121}.

\bibitem{gh-mo}
\refer{N. Ghoussoub, A. Moradifam:} {Bessel pairs and optimal Hardy  
and Hardy-Rellich inequalities,}{ Math. Ann.,} { 349} {(2011), 1--57}.




\bibitem{Kova-Pere-Seme}
\refer{V. F. Kovalenko, M.A. Perelmuter, Ya. A. Semenov:} {Schr\"odinger operators with $L_w^{l/2}(\R^l)$-potentials,} {J. Math. Phy.,} {22, n.5} {(1981),  1033-1044}.


\bibitem{met-pal-wack}
\refer{G. Metafune, D. Pallara, M. Wacker:} {Feller Semigroups on
$\R^N$,} {Semigroup Forum,} {153} {(2002),  179-206}.

\bibitem{met-spi}
\refer{G. Metafune, C. Spina: }{An integration by parts formula in Sobolev spaces, }{Mediterranean Journal of Mathematics, }{Vol. 5, N. 3, }{(2008), 359-371}.



\bibitem{met-spi2}
\refer{G. Metafune, C. Spina: }{Elliptic operators with unbounded coefficients in $L^p$ spaces, }{ Annali Scuola Normale Superiore di Pisa Cl. Sc. (5), }
{Vol XI}{ (2012), 2285-2299 }.

\bibitem{met-spi1}
\refer{G. Metafune, C. Spina: }{A degenerate elliptic operators with unbounded coefficients, }{ preprint (2012).}{ }{ }

\bibitem{met-spi-tac}
\refer{G. Metafune, C. Spina, C. Tacelli: }{ Elliptic operators with unbounded diffusion and drift
coefficients in $L^p$ spaces,}{ Advances Diff. Equat., }
{}{ (to appear)}.

\bibitem{mitidieri}
\refer{E. Mitidieri: }{A simple approach to Hardy inequalities, }{Mathematical Notes, }
{67 N. 4}{ (2000), 479-486}.


\bibitem{muck-wheed}
\refer{B. Muckenhoupt, R. Wheeden:} {Weighted Norm Inequalities for Singular and Fractional Integrals,} {Transactions of the America Mathematical Society,} {Vol. 161,} {(1971),  249-258}.

\bibitem{musina}
\refer{R. Musina:}{Optimal Rellich-Sobolev constants and their extremals,}{preprint (2013).}{}{}

\bibitem{musina1}
\refer{R. Musina:}{Weighted Sobolev spaces of radially symmetric functions,}{preprint (2013).}{}{}

\bibitem{nagel}
\refer{R. Nagel ({\rm ed}):} {One-parameter Semigroups of Linear Operators,}{Lecture Notes in Mathematics,}{Vol. 1184}{(1980), 1-24}.

\bibitem{okazawa}
\refer{N. Okazawa:} {$L^p$-theory of Schr\"odinger operators with strongly singular potentials,}{Japan. J. Math.,}{22}{(1996), 199-239}.


\bibitem{ou}
\refbook{E. M. Ouhabaz:} {Analysis of Heat Equations on Domains,}
{Princeton University Press}.


\bibitem{pazy}
\refbook{A. Pazy: }{Semigroups of linear operators and applications to partial differential equations, } {Applied mathematical sciences {\bf 44}, New York : Springer-Verlag, (1983)}.

\bibitem{rellich}
\refer{F. Rellich: } {Halbbeschr\"ankte Differentialoperatoren h\"oherer Ordnung,}{Proceedings of the International Congress of Mathematicians,}{Vol. III}{(1954), 243-250}.


\bibitem{stein}
\refer{E. M. Stein:}{Note on Singular integrals,}{Proc. Amer. Math. Soc.,}{8}{(1957), 250-254}.


\bibitem{stein-weiss}
\refbook{E. M. Stein, G. Weiss:} {Introduction to Fourier analysis on euclidean spaces,} {Princeton University Press ,} {(1971).}

\bibitem{vilenkin}
\refbook{N. Ja. Vilenkin:} {Fonctions sp\'{e}ciales et th\'{e}orie de la repr\'{e}sentation des groupes,} {Dunod Paris ,} {(1969).}

\bibitem{Voigt}
\refer{J. Voigt: }{One-parameter semigroups acting simultaneously on different $L_p$-spaces, }{ Bull Soc. Roy. Sci. Liege }
{61}{ (1992), 465-470 }.

\end{thebibliography}
\end{document}